\documentclass{ifacconf}

\usepackage{graphicx}
\usepackage{natbib}
\usepackage{amsmath,amsfonts,amsopn}
\usepackage{bbold}
\usepackage{amssymb}
\usepackage{nccmath}
\usepackage{color}
\usepackage{latexsym}
\usepackage{float}
\usepackage{mathrsfs} 
\usepackage{xcolor}
\usepackage{stmaryrd} 
\usepackage{tcolorbox}
\usepackage{mathrsfs}
\usepackage{enumitem}
\usepackage{multirow,multicol}
\usepackage{lineno}
\usepackage{verbatim,bm}
\usepackage{booktabs}
\usepackage{algorithm} 
\usepackage{algorithmic}
\usepackage[algo2e,vlined,linesnumbered,ruled]{algorithm2e}
\usepackage{mathtools}

\allowdisplaybreaks

\DeclarePairedDelimiter{\abs}{\lvert}{\rvert}

\theoremstyle{plain}
\newtheorem{theorem}{Theorem}[section]
\newtheorem{corollary}[theorem]{Corollary}
\newtheorem{proposition}[theorem]{Proposition}

\theoremstyle{definition}
\newtheorem{definition}[theorem]{Definition}
\newtheorem{remark}[theorem]{Remark}
\newtheorem{example}[theorem]{Example}

\newcommand{\C}{\mathbb{C}}
\newcommand{\ord}[2]{\operatorname{ord}_{#2}\left(#1\right)}

\begin{document}
\begin{frontmatter}

\title{Stability, Delays and Multiple Characteristic Roots in Dynamical Systems: A Guided Tour} 
\author[First]{Silviu-Iulian Niculescu} 
\author[First,Second]{Islam Boussaada} 
\author[Third]{Xu-Guang Li} 
\author[First]{Guilherme Mazanti}
\author[Fourth]{C\'{e}sar-Fernando M\'{e}ndez-Barrios} 

\address[First]{Universit\'e Paris-Saclay, CNRS, Inria, CentraleSup\'elec, Laboratoire des Signaux et Syst\`emes, 91190, Gif-sur-Yvette, France (e-mail: \{Silviu.Niculescu, Islam.Boussaada, Guilherme.Mazanti\}@l2s.centralesupelec.fr).}
\address[Second]{Institut Polytechnique des Sciences Avanc\'{e}es (IPSA), 63 boulevard de Brandebourg, 94200 Ivry-sur-Seine, France}
\address[Third]{School of Information Science and Engineering, Northeastern University, Shenyang, Liaoning 110004 China (e-mail: masdanlee@163.com).}
\address[Fourth]{Facultad de Ingenier\'{\i}a, Universidad Aut\'{o}noma de San Luis Potos\'{\i} (UASLP), San Luis Potos\'{\i}, M\'{e}xico (e-mail: cerfranfer@gmail.com)}

\begin{abstract}   
This paper presents a guided tour of some specific problems encountered in the stability analysis of linear dynamical systems including delays in their systems' representation. More precisely, we will address the characterization of multiple roots of the corresponding characteristic function with a particular emphasis on the way these roots are affected by the system's parameters and the way that they can be used to control. The paper covers several approaches (perturbation techniques, hypergeometric functions) leading to some methods and criteria (frequency-sweeping, multiplicity-induced-dominancy) that can be implemented (software toolboxes) for analyzing the qualitative and quantitative properties induced by the delays and other parameters on the system's dynamics. A particular attention will be paid to the so-called partial pole placement method based on the multiplicity-induced-dominancy property. The presentation is as simple as possible, focusing more on the main intuitive ideas and appropriate mathematical reasoning by analogy in the presentation of the theoretical results as well as their potential use in practical applications. Illustrative examples complete the paper. 
\end{abstract}


\end{frontmatter}

\section{Introduction}

Time-lag or dead-time, aftereffect, post action, deviating or lagging argument, heredity or hereditary effects represent some of the existing synonyms in the literature to describe the presence of a (time-)delay in the mathematical models of systems' dynamics. Large classes of physical, chemical and/or biological processes where the heterogeneity of the temporal phenomena needs a deeper understanding of the system's behavior make use of \emph{delays\/} in order to better capture the underlying mechanisms of such processes. Among others, one may cite: transport and propagation in interconnected cyber-physical systems subject or not to communication constraints or  incubation periods, maturation times, age structure, seasonal/diurnal variations in epidemic models. Such systems belong to the class of infinite-dimensional systems and there exist several ways to represent their dynamics. Functional differential equations (FDEs) sometimes called delay-differential equations (DDEs) are, by now, a classical framework for studying the qualitative and quantitative effects induced by the delays on the dynamics. Throughout the paper, we will adopt such a model representation. 
For a history of DDEs, in our opinion, the paper by \cite{hale:06} captures the most important results concerning the preservation of stability/instability of equilibria under small (nonlinear) perturbations in both frequency- (semigroups) and time-domain (mainly Lyapunov). Intuitively, extending, adapting ideas from ordinary differential equations (ODEs) to DDEs was, in some sense, the natural way to develop the corresponding theory and results concerning existence, uniqueness and continuation of solutions, the dependence on parameters\footnote{including the delays} are similar to the one derived for ODEs with a few additional technicalities\footnote{due to the infinite-dimensional character of the DDEs} as pointed out in \cite{Hale1993}. There exists an abundant literature on these topics in Mathematics, Physics, Engineering, Economics and Life Sciences since the 1970s (for older references, see, for instance, the excellent annotated bibliography by \cite{weiss:59} completed, a few months later, by \cite{chosky:60}, both published in Control area). 

In the sequel, we will focus on some problems arising in the stability analysis of linear DDEs in a frequency-domain framework. For a basic construction of the elementary solution of a linear DDEs using the Laplace transform as well as fundamental spectral properties of DDEs, we refer to \cite{kappel:06}. As briefly explained in Section \ref{sec:pre}, and similar to linear ODEs, the location of the roots of the characteristic function is essential to conclude on the (asymptotic) stability of linear DDEs. In particular, the way the system's parameters may affect the roots location is important in both quantitative and qualitative analysis of the system's behavior. 

In the 1940s, the construction of the parametric-plot of the frequency response of open-loop single-input/single-output (SISO) linear time-invariant (LTI) and the use of the Nyquist criterion to analyze the asymptotic stability analysis of  the corresponding closed-loop system became extremely popular, easy to understand and to apply and it is at the origin of most of the \emph{graphical (stability) approaches\/} and \emph{tests\/} in the open literature. Excepting the extension of such ideas to deal with  some particular classes of nonlinear systems (mainly, the well-known circle and Popov criteria), a particular attention was devoted to LTI systems including one delay in the input/output channel, and a lot of results published in delay area in the 50s--70s deal with various extensions of Nyquist and Michailov criteria (see, for instance, \cite{krall:65}). For a deeper discussion of the Michailov criterion applied to delay systems, we refer to \cite{barker:79} and for some connections between Michailov and Nyquist criteria applied to delay systems we mention the almost forgotten paper by \cite{chen-tsay:76}. Further graphical tests include the well-known \emph{root-locus methods}\footnote{whose origins go back to the works of Evans at the end of the 1950s see, e.g. \cite{evans:50} and the references therein} and the Satche's diagrams (see, e.g. \cite{satche:49}), where the last ones can be interpreted as a variant of the Nyquist criterion and are sometimes called dual root-locus methods. 
For a pedagogical presentation as well as some extensions of the root locus methods and Nyquist criterion to deal with SISO LTI systems with one delay in the input/output channel, we refer to \cite{krall:68}\footnote{see also \cite{krall:70} for a survey of the root-locus methods}. Finally, to complete the discussion above, it is worth mentioning the extension of Popov criteria to delay systems (see, for instance, \cite{halanay:66}, \cite{rasvan:75}).
 
At the beginning of the 80s, \cite{els:73} mention three tests for checking the asymptotic stability of DDEs: the amplitude-phase method (referring to Tsypkin's contributions), the so-called $D$-partition method and the direct generalization of the Routh--Hurwitz method (mainly \u{C}ebotarev's contributions), with a deeper discussion of the first two methods by using the argument principle as well as the Rouch\'{e}'s lemma.

The $D$-partition method has a long history started with the contribution of Vishnegradsky in the 1880s and its application to DDEs goes back to the 1950s (see, e.g. \cite{Neimark}). For a deeper discussion on the construction of stability charts for low-order systems including delays in their system representation, we refer to \cite{Pinney:58} and \cite{Stepan89}. The so-called amplitude-phase method, proposed by \cite{Tsypkin:46}, is at the origin at most of the existing \emph{frequency-sweeping tests\/} in the open literature and it will be briefly addressed in Section \ref{sec:pre}.  However, to the best of the authors' knowledge, the notion of ``frequency-sweeping'' was first mentioned \cite{cl:95} and \cite{c:95} into a different methodological frame. Excepting these methods, at the end of the 70s, \cite{lee-hsu:69} proposed the so-called \emph{$\tau$-partition\/} method that can be seen as the ``dual'' of the $D$-partition method if one takes into account the way the coefficients and the delay(s) are treated\footnote{$D$-partition: fixed delays with all the other parameters free and $\tau$-partition: fixed coefficients and the delays are free parameters}. For its extension in the most general case for analytical functions with respect to one delay parameter, we refer to \cite{cd:86}.  Independently, a similar analysis was proposed by \cite{Walton-Marshall} one year later, and both contributions are at the origin of various results published in the literature in the last 20 years.

The monograph by \cite{cm:49} collects some of the developments on the Routh--Hurwitz problem in the period 30s--50s, including their direct generalization to quasipolynomials. Since most of these methods are difficult to apply, starting with the 90s, the use of some appropriate bilinear transformations allowed to reformulate and to simplify the detection of characteristic roots on the imaginary axis. Such ideas are at the origin of the so-called \emph{pseudodelay techniques} initiated by \cite{reka}, that have been further refined by \cite{2tho0} and, more recently, by \cite{olgac} leading to some new methods.

The development of robust control methodologies in the 90s allowed to reconsider some of the ideas above in the new frame and to develop some new techniques. In particular, the interpretation of delays as uncertain parameters is at the origin of a series of contributions on the so-called \emph{delay-independent/delay-dependent\/} stability\footnote{To the best of the authors' knowledge, the concepts of stability/stabilization ``independent of delay'' were formally introduced by \cite{Kamen:82}.} problem with a particular attention to the computation of the \emph{delay margin\/} in the delay-dependent case\footnote{To the best of the authors' knowledge, the notion of \emph{delay margin\/} was introduced by \cite{c:95} and \cite{cgn:95}.}. For further discussions on such topics as well as a long list of references, we refer to \cite{n:01-springer},  \cite{keqin-book} and \cite{fridman:14} (see also \cite{csm:11-delay} and \cite{richard:03}).

Since the earlier studies devoted to the analysis of linear DDEs in the 50s, the existence of multiple characteristic roots was mentioned by \cite{Hayes50} (scalar DDE) \cite{Pinney:58} (scalar and second-order DDEs) in characterizing the stability charts in the corresponding parameter space, but without any attempt to develop an appropriate methodology in the general case. Extending the ideas of \cite{Stepan79} (see also \cite{Stepan89}), \cite{Hassard1997} proposed an analytical criterion to count the unstable roots. Such a criterion, based on the argument principle, takes into account multiple roots on the imaginary axis subject to some appropriate constraints. 
Next, as mentioned in \cite{Hale1993}, deriving an abstract perturbation theory for DDEs is not a trivial problem. However, by exploiting the structure of the characteristic functions, \cite{cfnz:10-siam-II} (see also the companion paper \cite{cfnz:10-siam-I} and \cite{CM:17}) derived an \emph{eigenvalue perturbation approach\/} for characterizing the asymptotic behavior of multiple characteristic roots located on the imaginary axis for some classes of DDEs\footnote{retarded DDEs with multiple commensurate delays}, and the underlying ideas are at the origin of some of the developments in the literature during the last years (see also the discussions in \cite{mn14:siam}). Finally, in the case of complete regular splitting, we mention \cite{SIMAX-MBN-17} (see also \cite{MM:19} for some insights in the classification).
  
The contributions of the paper are twofold: first, to offer a guided tour on some of the approaches developed by the authors during the last decade on the analysis of the effects induced by multiple characteristic roots on the system's dynamics. Second, we are interested to emphasize the way that such findings can be used in control. To improve the readability of the paper, most of the presented results are completed by illustrative examples. 

The paper is organized as follows: the problem statement and some prerequisites are briefly presented in Section \ref{sec:prob} and two motivating examples in Section \ref{sec:motiv}. Sections \ref{sec:Pui}--\ref{sec:mult-del} focus on some approaches to characterize the multiple characteristic roots located on the imaginary axis and are strongly related to the so-called $\tau$-partition methods. Next, Section \ref{sec:MID} covers an interesting property valid only for delay systems --- \emph{multiplicity-induced-dominancy\/} --- opening interesting perspectives in control. Finally, some concluding remarks end the paper.

\emph{Notations\/}: In this paper, $\mathbb{N}^\ast$ denotes the set of positive integers and $\mathbb{N} = \mathbb{N}^\ast \cup \{0\}$. The set of all integers is denoted by $\mathbb{Z}$ and, for $a, b \in \mathbb{R}$, we denote $\llbracket a, b\rrbracket = [a, b] \cap \mathbb{Z}$, with the convention that $[a, b] = \emptyset$ if $a > b$. For a complex number $\lambda$, $\Re(\lambda)$ and $\Im(\lambda)$ denote its real and imaginary parts, respectively; $\mathbb{C}_-$ and  $\mathbb{C}_+$ denote the sets $\{ \lambda  \in
\mathbb{C}: \Re(\lambda ) < 0\} $ and $\{ \lambda \in \mathbb{C}: \Re (\lambda ) > 0\} $, respectively; furthermore, $j\mathbb{R}$ (with $j=\sqrt{-1}$) denotes the imaginary axis. The \textit{order} 
of a power series $f(x,y)=\sum_{i,k} a_{i,k}x^{i}y^{k}$ will be denoted by $\ord{f}{}$ and defined as the smallest number $n=i+k$ such that $a_{i,k}\neq0$. The order of the power series $f$ with respect to the variable $x$, will be denoted by $\ord{f}{x}$ and defined similarly. For $n\in\mathbb{N}$, let $\overrightarrow{x}:=\left(x_1,x_2,\dotsc,x_n\right)$, then the ring of complex formal power series will be denoted by $\C[[\overrightarrow{x}]]$, with subring $\C\{\overrightarrow{x}\}$ of convergent power series. 

\section{Parameters, delays and dynamics: Problem formulation and prerequisites}
\label{sec:prob}

Consider a positive integer $n_p\in \mathbb{Z}_+$, and an open set $\mathcal{O}\in\mathbb{R}^{n_p}$. For a set of parameters $\overrightarrow{p}\in\mathcal{O}$, introduce the transcedental complex-valued function $\Delta:\mathbb{C}\times \mathcal{O}_p \mapsto \mathbb{C}$ defined by:
\begin{equation}
\label{quasi-gen}
\Delta(\lambda; \overrightarrow{p},\overrightarrow{\tau}) := P_0(\lambda,\overrightarrow{p}) +\sum_{i=1}^{n_d}P_i(\lambda,\overrightarrow{p})e^{-\lambda \tau_i(\overrightarrow{p})},
\end{equation}
where the components $\tau_i:\mathcal{O} \mapsto \mathbb{R}_+$, $i\in \llbracket 1, n_d\rrbracket$ of the delay vector $\overrightarrow{\tau}\in\mathbb{R}^{n_d}$ are assumed to be sufficiently smooth, non-negative and bounded functions for all the parameters $ \overrightarrow{p}\in\mathcal{O}$. Next, $P_i$, ($i\in \llbracket 0, n_d\rrbracket$) denote real polynomials in the complex variable $\lambda$ depending on the parameters $ \overrightarrow{p}$, such that the applications $\overrightarrow{p}\mapsto P_i(\cdot,\overrightarrow{p})$ are well-defined and sufficiently smooth on the open $\mathcal{O}$. Assume further that $P_0$ is a monic\footnote{The leading coefficient $p_{0,n}$ of $P_0$ is equal to one: $p_{0,n}=1$.} polynomial and that $n=\mathrm{deg}(P_0)>\mathrm{max}_{\overrightarrow p}\mathrm{deg}(P_i)$ for all $\overrightarrow{p}\in\mathcal{O}$ and $i\in \llbracket 1, n_d\rrbracket$ with the observation that we are not excluding the cases when leading coefficients of the polynomial $P_i$ may vanish for some values of the parameters $\overrightarrow{p}$. These last assumptions allow guaranteeing that the quasipolynomial (if it exists?!) is always of \emph{retarded type}.

The transcendental function $\Delta$ in (\ref{quasi-gen}) covers a lot of cases encountered in the analysis of delay systems depending on parameters. Such parameters may define some particular structure of the coefficients of the polynomials $P_i$ ($i\in \llbracket 0, n_d\rrbracket$) and/or of the delays\footnote{for instance, the case of rationally-dependent delays}, some dependence between the coefficients of the polynomials and the delays, or it may reflect the way the controller's gains appear in the characteristic function of the closed-loop system. 

To fix better the ideas, consider the strictly proper LTI SISO  linear system $\Sigma(A,b,c^T)$ with the state-space representation:
\begin{equation}
\label{siso}
\Sigma: \quad \left\{ 
\begin{aligned}
& \dot x(t)  = Ax(t)+ b u(t)\\
& y(t) = c^Tx(t),
\end{aligned} 
\right.
\end{equation}
where the transfer function $H_{yu}(\lambda)$ of $\Sigma$ writes as $H_{yu}(\lambda)=P_1(\lambda)/ P_0(\lambda)$, for some appropriate real polynomials $P_i$, $i \in \llbracket 0,1 \rrbracket$, whose coefficients are given by the ``entries'' of $\Sigma$. Assume now that $\Sigma$ is controlled by the delayed output feedback $u(t)=-ky(t-\tau)$ with $k\in\mathcal{O}_p\subset \mathbb{R}^*$. In our case, the pair $(k,\tau)$ simply represents a ``delay-block'' (controller). Under the assumption that the pair $(k,\tau)\in \mathcal{O}_p\times \mathbb{R}_+$ defines our \emph{parameters\/} then, some simple computations show that the stability of the system in closed-loop reduces to the analysis of the location of the spectrum of the quasipolynomial $\Delta(\cdot; k,\tau)$ given by:
\[
\Delta(\lambda; k,\tau):=P_0(\lambda) +kP_1(\lambda)e^{-\lambda\tau},
\]
In the case when the gain verifies $\abs{k} = 1$, we refer to \cite{Tsypkin:46} for a first discussion on the stability with respect to the delay parameter. Finally, it is easy to observe that we arrive to the same characteristic quasipolynomial $\Delta$ for the closed-loop system if the delay $\tau$ is in the input channel and not in the output.
 
Using the same terminology as \cite{BC63}, the zeros of $\Delta$ are called \emph{characteristic roots}. 
Denote by $\sigma_{s}(\Delta)\subset \mathbb{C}$ the whole set of characteristic roots of $\Delta$.
With these notations and notions, the stability problem can be formulated as follows: \emph{find the whole set of parameters $\overrightarrow{p}\in\mathcal{O}$ guaranteeing that the corresponding characteristic roots are located in $\mathbb{C}_-$.} 

\subsection{Characteristic roots: prerequisites}
\label{sec:pre}

Under the assumption that $P_0$ is a monic polynomial, we will further assume that the delay  $\overrightarrow{\tau}$ and the parameter $\overrightarrow{p}$ vectors are \emph{independent\/} each-other. Thus, in the retarded case, for a set of parameters $(\overrightarrow{p},\overrightarrow{\tau})\in\mathcal{O}_p \times  \mathbb{R}_+^{n_d}$,  the characteristic function $\Delta:\mathbb{C}\times \mathcal{O}_p \times 
\mathbb{R}_+^{n_d} \mapsto \mathbb{C}$ writes as:
 \begin{equation}
\label{quasi-gen-param}
\Delta(\lambda; \overrightarrow{p},\overrightarrow{\tau}) := P_0(\lambda,\overrightarrow{p}) +\sum_{i=1}^{n_d} P_i(\lambda,\overrightarrow{p})e^{-\lambda \tau_i},
\end{equation}
where  $\deg(P_i)<\deg(P_0)$, for all $i \in \llbracket 1,n_d \rrbracket$ and for all $\overrightarrow{p}\in \mathcal{O}_p$. 
A particular important case of (\ref{quasi-gen-param}) is represented by the class of quasipolynomials with \emph{commensurate delays}, that can be simply written by taking $\tau_i=i\tau$, for all $i\in  \llbracket 1,n_d \rrbracket$. In the simplest case when $n_d=1$, (\ref{quasi-gen-param}) rewrites as:
\begin{equation}
\label{quasi-1}
\Delta(\lambda; \overrightarrow{p},\tau) := P_0(\lambda,\overrightarrow{p}) +P_1(\lambda,\overrightarrow{p})e^{-\lambda \tau}.\end{equation}
The quasipolynomial $\Delta$ given by (\ref{quasi-1}) has some nice and interesting properties that will be exploited in the sequel. For instance, any vertical stripe of the complex plain includes a finite number of characteristic roots. Furthermore, there exists a real number $\gamma$, such that all the characteristic roots are confined to the half-plane $\mathbb{C}_\gamma$: $\left\{\lambda\in \mathbb{C}: \Re(\lambda)<\gamma \right\}$ (see, for instance, \cite{mn14:siam} and the references therein).

Based on Rouch\'{e}'s lemma (see, e.g., \cite{Ahlfors:79}), we have the following:

\begin{theorem}
\label{theo:cont}
Under the assumption that $P_0$ is monic and $\deg(P_0)>\deg(P_1)$ for all $\overrightarrow{p}\in \mathcal{O}_p$, let $\lambda_0$ be a characteristic root of the quasipolynomial $\Delta(\cdot; \overrightarrow{p_0}, \tau_0)$ with multiplicity $k$. Then there exists a constant $\bar\varepsilon>0$ such that for all $\varepsilon>0$
satisfying $\varepsilon<\bar\varepsilon$, there exists a $\delta_\varepsilon>0$ such that 
$\Delta(\lambda;\ \overrightarrow{p_0}+\delta \overrightarrow{p_0}, \tau_0+\delta\tau_0)$, where $\delta\tau_0\in\mathbb{R}$, $\abs{\delta\tau_0} < \delta_\varepsilon$, $\tau_0+\delta\tau_0\geq 0$, 
 $\delta \overrightarrow{p_0}\in\mathbb{R}^{n_p}$, $\|\delta \overrightarrow{p_0}\|_2<\delta_\varepsilon$, $\overrightarrow{p_0}+\delta \overrightarrow{p_0}\in\mathcal{O}_p$
has exactly $k$ zeros (multiplicity taken into account) in the disc
$\left\{\lambda\in\mathbb{C}: \ |\lambda-\lambda_0|<\varepsilon\right\}$.  \end{theorem}

\begin{remark}[``Small'' delays case]
\label{rem:small-del}
Consider now the case $\Delta(\lambda;\tau)=P_0(\lambda)+P_1(\lambda)e^{-\lambda\tau}$ with $\mathrm{deg}(P_0)>\mathrm{deg}(P_1)$, for a sufficiently small delay value $\tau=\varepsilon>0$. The finite roots of $\Delta(\cdot;\varepsilon)$ can be made arbitrarily close to the finite roots of $\Delta(\cdot;0)$ and there exists an infinite number of  roots whose real parts approach negatively infinite\footnote{For an elementary proof, see, e.g., the Appendix of \cite{sha-ka:69}.}. In other words, when increasing the delay from $0$ to $0_+$, although the system changes its character\footnote{from finite- to infinite-dimensional}, the stability/instability of the delay-free system is preserved for sufficiently small delays. However, such a property does not necessarily hold in all the cases, and there are two particular situations of interest:
(i) \emph{neutral\/} case\footnote{not addressed in this paper} ($\mathrm{deg}(P_0)=\mathrm{deg}(P_1)$) and (ii)  \emph{delay-dependent coefficients\/} of $P_0$, $P_1$.\\
The last case may appear in the PD-control of LTI SISO systems when the derivative action is implemented by using an Euler delay-difference approximation scheme, see, e.g. \cite{mnmr:22}. 
For further discussions on delay systems including delay-dependent coefficients, we refer to \cite{Chi:18-over} and the references therein. $\Box$
\end{remark}

\subsection{Spectral abscissa function: definition and properties}
\label{sub-sec:spectral-abs}

For the analysis of stability, it is important to know where the rightmost characteristic root is located as well as the way it is affected by parameters change. To answer to such questions, introduce now the \emph{spectral abscissa function\/} $(\overrightarrow{p},\tau)\in\mathcal{O}_p\times \mathbb{R}_+ \mapsto \alpha_s(\overrightarrow{p},\tau)\in\mathbb{R}$ defined by
\
\[
\alpha_s(\overrightarrow{p},\tau):= \sup\left\{\Re(\lambda):\ \ \Delta(\lambda;\ \overrightarrow{p},\tau)=0, \overrightarrow{p} \in \mathcal{O}_p \right\}.
\]
As a consequence of Theorem \ref{theo:cont}, we have two properties:
\begin{itemize}
\item[(i)] If $\mathrm{deg}(P_0)>\mathrm{deg}(P_1)$, $\alpha_s$ always \emph{exists}, is \emph{bounded} and \emph{continuous}.
\item[(ii)] As the delay and/or parameters vary, the multiplicity summation of the roots of $\Delta$ in open $\mathbb{C}_+$ can change only \emph{if a root appears on} or \emph{crosses} the \emph{imaginary axis}\footnote{For an elementary proof of such a property, we refer to \cite{cg:82} in the case of a second-order system with respect to the delay parameter}. 
\end{itemize}

Thus, understanding the behavior of the characteristic roots located on $j\mathbb{R}$ with respect to the parameters' change becomes essential for a complete characterization of the stability regions in the corresponding parameter-space.

\begin{remark}
The ideas above still hold in the commensurate delays case ($\tau_i=i\tau$ for $i\in \llbracket 1,n_d \rrbracket$ and $\tau \in \mathbb{R}_+$ in (\ref{quasi-gen-param})). For incommensurate delays, by introducing \emph{delay rays\/} $\left\{r\overrightarrow{\tau}: r\in \mathbb{R}_+\right\}$, \cite{datko:78} proved that the continuity of the spectral abscissa holds with respect to $r\in\mathbb{R}_+$. $\Box$
\end{remark}

\begin{remark}
To construct the stability charts in the scalar and second-orded DDEs, \cite{Pinney:58} introduced the so-called \emph{$(x_r,k_r)$ plateau} set, that is the set of parameters for which $\Delta$ has $k_r$ and only $k_r$ roots of real part greater than $x_r$. Thus, $(0,0)$-root plateau corresponds to the \emph{stability regions}, and the minimal value of $x_r$ of the $(x_r,0)$-root plateau corresponds to the \emph{spectral abscissa}. $\Box$
\end{remark}

\subsection{Hyperbolicity, switches and reversals}

Reconsider the SISO system (\ref{siso}) in closed-loop under the assumption that $k=1$, i.e., $\Delta$ rewrites as: $\Delta(\lambda;\tau):=P_0(\lambda)+P_1(\lambda)e^{-\lambda\tau}$ and let us focus on the $\tau$-partition. Assume further that $P_0$ and $P_1$ are coprime. If 
 \begin{equation}
 \label{FS-cond-hyp}
     \abs{P_1(j\omega)} < \abs{P_0(j\omega)}, 
 \end{equation}
for all $\omega \in \mathbb{R}$, then $\sigma_s(\Delta)\cap j\mathbb{R} =\emptyset$.
By using Theorem~\ref{theo:cont}, it follows that the characteristic roots of $\Delta$ can not migrate from $\mathbb{C}_-$ to $\mathbb{C}_+$ or vice-versa if $\tau$ is increased from $0$ to $+\infty$. Such a system is called \emph{hyperbolic} and it has an interesting property: the location of the spectrum of the polynomial $P_0+P_1$ will define the stability/instability of the system for all delays $\tau\in\mathbb{R}_+$.  For further discussions in a more general setting, we refer to \cite{hit:85}\footnote{For the characterization of the commensurate delays case, see, e.g., \cite{n:01-springer}.}.
 
Next, it is easy to see that if $0\in \sigma_s(\Delta(\cdot;0))$, then $\Delta(0;\tau)=0$, $\forall\tau\in\mathbb{R}_+$. Thus, the origin will be an \emph{invariant root}\footnote{The common roots $P_0$ and $P_1$ on $j\mathbb{R}$ are also invariant roots with respect to $\tau$.}. Now, if $0\not\in\sigma(P_0+P_1)$, checking (\ref{FS-cond-hyp}) for $\forall\omega\in\mathbb{R}_+^*$ is sufficient to guarantee hyperbolicity. 
Assume now that $\sigma(P_0+P_1)\subset\mathbb{C}_-$.
As observed by \cite{Tsypkin:46}, the closed-loop system is \emph{delay-independently stable} if and only if the condition (\ref{FS-cond-hyp}) holds for all $\omega\in\mathbb{R}_+^*$, and it can be simply checked from the plot of $z_1$, where the application $\omega \mapsto z_1(\omega):=-P_0(j\omega)/P_1(j\omega)$, for $\omega\in\mathbb{R}_+^*$ defines the simplest \emph{frequency sweeping curve (FSC)}.

Consider now the case when the closed-loop system is not hyperbolic. Then there exists at least one value $\omega_c\in\mathbb{R}$, such that $\Delta(j\omega_c;\tau)=0$ for some delay $\tau=\tau_c\in\mathbb{R}_+$. Such a frequency $\omega_c$ will be called \emph{crossing frequency}, and the collection of all ``$\omega_c$'' will define the \emph{crossing set\/}:
\begin{equation}
\label{xing-set}
\Omega_c := \left\{ \omega \in \mathbb{R}: \; \abs{P_0(j\omega)} = \abs{P_1(j\omega)} \right\}.    
\end{equation}
At this stage, there are two important remarks:
\begin{itemize}
    \item[(i)] first, $\mathrm{card}(\Omega_c)$ is \emph{finite}, and its computation reduces to the computation of the positive roots of an appropriate polynomial;
    \item[(ii)] second, the knowledge of a crossing frequency $\omega_{i,c} \in \Omega_c$ will allow to compute the minimal critical delay value $\tau_{i,c}^*\in\mathbb{R}_+$\footnote{Such a value always exists and it may be $0$.} that will generate the set of \emph{(critical) crossing delays}
\begin{equation}
\label{del-set}
    \mathcal{T}(\omega_{i,c}):=\left\{ \tau_{i,c}^*+\frac{2k\pi}{\omega_c}\geq 0, \quad k \in \mathbb{Z}\right\}.
\end{equation}
\end{itemize}
For a deeper discussion of the remarks (i)--(ii) above, see, for instance, \cite{mn14:siam}.

Under the assumption of a simple characteristic root $\omega_0\in\Omega_c$ for some delay $\tau_0\in \mathcal{T}(\omega_0)$, \cite{cg:82} discussed the behavior of the characteristic root $j\omega_0$ for values close to $\tau_0$ by using the ``quantity'' $\mathrm{sgn}(\Re(d\lambda)d\tau)$ evaluated a $\lambda=j\omega_0$ and $\tau=\tau_0$. Such an idea was further refined in \cite{cd:86} and largely use in the open literature during the last 30 years. More precisely, if the characteristic root located on the imaginary axis moves towards instability (stability), we will have a \emph{stability switch\/} (\emph{reversal\/})\footnote{To the best of the authors' knowledge, during the 80s, the notions of \emph{(stability) switches/reversals\/} appear in Cooke's publications.}. Finally, the degenerate case when $\Re(d\lambda)d\tau=0$ needs to take into account higher-order derivatives (see, for instance, \cite{sha-ka:69} and the references therein).

\subsection{Quasipolynomial degree and multiplicity}
\label{sub-sec:Van}

Recall the general quasipolynomial (\ref{quasi-gen-param}). The integer $\mathscr{D}_{PS} = n_d + \sum_{k=0}^{n_d} m_k$ is called the \emph{degree} of $\Delta$ (see, for instance, \cite{Wielonsky2001Rolle}).

\begin{remark}
\label{remk:PolyaSzego}
A classical result known as \emph{P\'{o}lya--Szeg\H{o} bound}, see, e.g., \cite{PS72}\footnote{This result  was first introduced and claimed in the problems collection published in 1925 by G.~P\'{o}lya and G.~Szeg\H{o}. In the fourth edition of their book \cite[Part Three, Problem~206.2]{PS72}, G.~P\'{o}lya and G.~Szeg\H{o} emphasized that the proof was obtained in the meantime by N.~Obreschkoff using the principle argument, see \cite{Obreschkoff}.} allows to establish a direct link between the degree of a quasipolynomial and the number of its roots in horizontal strips of the complex plane.
As an immediate consequence, given a root $\lambda_0 \in \mathbb C$ of a quasipolynomial \eqref{quasi-gen-param}  of degree $\mathscr{D}_{PS}$, by letting the horizontal strip a line,  one concludes that any root of a quasipolynomial  has multiplicity at most $\mathscr{D}_{PS}$. $\Box$
\end{remark}

\begin{remark}
Using a constructive algebraic approach ba\-sed on functional Birkhoff matrices, \cite{BN-ACAP-2016} showed that the maximal admissible multiplicity of quasipolynomial's roots is the P\'{o}lya--Szeg\H{o} bound. Furthermore, in the lacunary case\footnote{when some coefficients of the quasipolynomial are identically zero}, it has been shown in \cite{BN-ACAP-2016} that the P\'{o}lya--Szeg\H{o} bound  cannot be reached and some sharper bounds for the admissible multiplicities has been established in some configurations. $\Box$
\end{remark}

\begin{remark}
The problem of identifying the maximal dimension of the eigenspace associated to a multiple singularity $\lambda=j\omega_0$ (with non-vanishing frequency $\omega_0\neq 0$) for time-delay systems as well as the explicit conditions guaranteeing such a configuration has been addressed in \cite{BN-TAC-2016}, and the conclusion is that the P\'{o}lya--Szeg\H{o} bound for the maximal admissible multiplicity is never reached when the crossing frequency is different from zero. $\Box$
\end{remark}

\section{Motivating examples}
\label{sec:motiv}

\subsection{Scalar case: Double zero singularity}

Consider the following scalar delay-differential system:
\begin{equation}\label{eq:scalar-aux-dde}
\dot y(t)+\alpha \left( y(t)-y(t-1) \right) =0,
\end{equation}
under appropriate initial conditions, where $\alpha\in \mathbb{R}^*$. The corresponding characteristic function rewrites as:
\begin{equation}\label{eq:scalar-aux-char-f}
\Delta (\lambda; \alpha) := \lambda+\alpha \left( 1 - e^{-\lambda} \right).
\end{equation}
It is easy to see that $\Delta (0; \alpha)=0$ for all $\alpha \in \mathbb{R}$, showing that such a root is \emph{invariant} with respect to $\alpha$. Furthermore, since $\Delta' (\lambda)=1+\alpha e^{-\lambda}$, then the root at the origin is \emph{double} if $\alpha=-1$. 
With these observations in mind, we have the following result:

\begin{proposition}
The scalar system (\ref{eq:scalar-aux-dde}) is unstable for all $\alpha \in (-\infty,-1)$, and the characteristic function $\Delta_\alpha$ has one strictly positive real characteristic root. If  $\alpha \in (-1,+\infty)$, excepting the root at the origin, the remaining characteristic roots of $\Delta$ (if any?!) are located in $\mathbb{C}_-$.
\end{proposition}

We have three important observations:
\begin{itemize}
\item First, when the real parameter $\alpha$ is increased from $-\infty$, one real characteristic root arrives from $+\infty$ and it will move on the real axis towards to $-\infty$ when $\alpha$ tends to $+\infty$. In other words, this root is \emph{``locked''\/} on the real axis for all values of the parameter $\alpha$ and the characteristic function $\Delta$ has always two real roots.
\item  Second, when $\alpha=-1$,  the system has a double characteristic root $\lambda=0$. It is easy to show that for general scalar DDEs including a single delay, the maximal multiplicity of a characteristic root is two and it can be reached only on the real axis. 
\item Finally, surprisingly, the double root at the origin is dominant in the sense that all the other roots are located in $\mathbb{C}_-$. 
 Such a property, called \emph{multiplicity-induced-dominancy\/} will be further addressed in the forthcoming sections. 
\end{itemize}

\begin{figure}[!ht]
\centering
\resizebox{0.9\columnwidth}{!}{\input{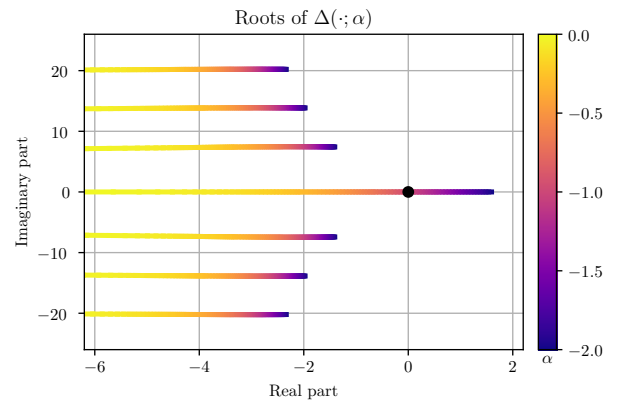}}
\caption{Roots of $\Delta(\cdot; \alpha)$ from \eqref{eq:scalar-aux-char-f} for $\alpha \in [-2, 0]$.}\label{fig:order-1}
\end{figure}

The behavior of the the roots of $\Delta(\cdot; \alpha)$ as $\alpha$ varies in the interval $[-2, 0]$ is illustrated in Fig.~\ref{fig:order-1}, in which different values of $\alpha$ in $[-2, 0]$ are represented by different colors and the root at $0$ for all $\alpha$ is represented by a black dot. We remark that, in addition to the root at $\lambda = 0$, $\Delta(\cdot; \alpha)$ has another distinct real root for $\alpha \in (-\infty, -1) \cup (-1, 0)$, which is positive if $\alpha < -1$ and negative if $\alpha > -1$. Notice that $\lambda = 0$ is the unique root in the case $\alpha = 0$ and, as $\alpha \to 0$, the real parts of all other roots of $\Delta(\cdot; \alpha)$ converge to $-\infty$, as previously described in Remark~\ref{rem:small-del}.

\begin{remark}
A deeper analysis of the existence of double roots in the scalar case can be found in \cite{nuss:02}. For further discussions the existence of real roots of the characteristic function for general scalar DDEs with respect to the system's coefficients we refer to a series of papers written by Wright at the end of the 60s (see, for instance, \cite{Wright1961Stability} and the references therein). Finally, to the best of the authors' knowledge, the first complete characterization of the stability regions in the scalar case (covering both retarded and neutral cases) can be found in \cite{Hayes50}. $\Box$
\end{remark}

\begin{remark}
The explicit expressions of the characteristic roots can be done by using the so-called \emph{Lambert W function\/} (see, for instance, \cite{corless:96} and the references therein), that is the (multivalued) inverse of complex function $\xi\in\mathbb{C}  \mapsto \xi e^{\xi} \in \mathbb{C}$. It has an infinite, but countable number of branches $W_k(\xi)\in \{w \in \mathbb{C}: \xi=we^w\}$, for $k \in \mathbb{Z}$. More precisely, the characteristic roots of $\Delta$ given by (\ref{eq:scalar-aux-char-f}) are expressed as:
$$
\lambda_k = -\alpha+ W_k\left(\alpha e^\alpha\right), \quad \forall k \in \mathbb{Z}.
$$
Each of these branches is locally analytic excepting the principal branch $W_0$ that is not differentiable at the point $\xi=-e^{-1}$, that corresponds to the case when the parameter $\alpha=-1$. A deeper discussion of the general scalar case by using the Lambert W function can be found in \cite{agu:03} (see also \cite{yng:10} for some extensions of these ideas to the analysis and synthesis of delay systems). $\Box$
\end{remark} 

\begin{figure}[!ht]
\centering
\resizebox{0.9\columnwidth}{!}{\input{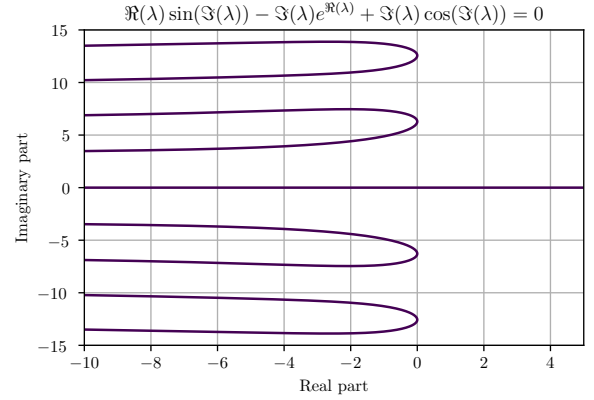}}
\caption{Roots in $\mathbb C$ of equation \eqref{eq:scalar-root-locus}.}\label{fig:root-locus}
\end{figure}

\begin{remark}
\label{remk:root-locus}
Observe that any root $\lambda$ $\Delta(\cdot; \alpha)$ necessarily satisfies:
\begin{equation}
\label{eq:scalar-root-locus}
\Re(\lambda) \sin(\Im(\lambda)) - \Im(\lambda) e^{\Re(\lambda)} + \Im(\lambda) \cos(\Im(\lambda)) = 0,
\end{equation}
and, conversely, if $\lambda \in \mathbb C$ satisfies \eqref{eq:scalar-root-locus} and $\Re(\lambda) \neq 0$, then there exists a unique $\alpha \in \mathbb R$ such that $\Delta(\lambda; \alpha) = 0$\footnote{The case of points on the imaginary axis satisfying \eqref{eq:scalar-root-locus} can be obtained as a limit as $\alpha \to \pm \infty$}. Fig.~\ref{fig:root-locus} represents the set of points in the complex plane satisfying \eqref{eq:scalar-root-locus}. Notice that the roots represented in Fig.~\ref{fig:order-1} all lie in the set of points represented in Fig.~\ref{fig:root-locus}. $\Box$
\end{remark}

\subsection{Inverted pendulum stabilization: Triple zero singularity}

Consider now a dynamical system modeling a friction free inverted pendulum on cart. The adopted model is studied in \cite{Krauskopf2004}, \cite{S2005}, \cite{BMN2015-SCL} and, in the sequel,  we keep the same notations.
\begin{figure}[!ht]
\begin{center}
\includegraphics[clip=true,trim=5mm 5mm 20mm 5mm,scale=0.9]{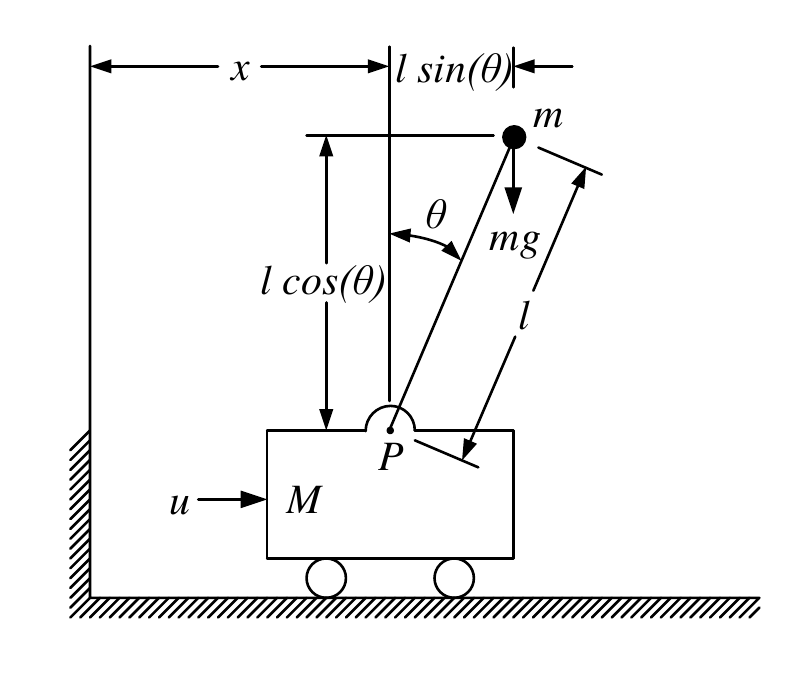}
\caption{Inverted pendulum on a cart.}\label{IP}
\end{center}
\end{figure}
In the dimensionless form, the dynamics of the inverted pendulum on a cart in Fig.~\ref{IP} is governed by the following second-order differential equation:
\begin{equation}\label{InvertedPendulum}
\!\!\left(\!\!1-\!\frac{3\epsilon}{4}\cos^{2}(\theta)\!\!\right)\ddot\theta+\frac{3\epsilon}{8}\dot\theta^2\sin(2\theta)-\sin(\theta)+u\cos(\theta)=0,
\end{equation}
where $\epsilon={m}/{(m+M)}$, $M$ the mass of the cart and $m$ the mass of the pendulum and $u$ represents the control law that is the horizontal driving force.
Consider now that such a system is controlled by using a standard delayed PD controller of the form $u(t)=k_p\,\theta(t-\tau)+k_d\,\dot\theta(t-\tau)$, with $(k_p,k_d)\in\mathbb{R}^2$. The local stability of the closed-loop system is reduced to study the location of the spectrum of the quasipolynomial $\Delta(\lambda; k_p,k_d,\tau):=Q(\lambda)+P(\lambda; k_p,k_d)e^{-\lambda\tau}$ where the polynomial $P$ is first-order and includes the gain parameters $(k_p,k_d)$. A generalized Bogdanov--Takens singularity with codimension three is identified in \cite{Krauskopf2004}. 

Consider now a simpler planar inverted pendulum in absence of friction:
\begin{equation}
 \label{pend-atay}
 \ddot \theta(t) -\frac{g}{l} \sin(\theta(t)) = u(t),
\end{equation}
where $\theta$ is the angular displacement\footnote{measured from the natural equilibrium position}, $g$ the gravitational acceleration, $l$ the pendulum length and $u$ the external torque. Assume that the controller\footnote{information available on the ``past'' (angular) position and not on the speed} includes two ``delay blocks'' $(k_i,\tau_i)$, with $i=1,2$, and the control law has the form: $u(t)=-k_1\theta(t-\tau_1)-k_2\theta(t-\tau_2)$. The characteristic function of the linearized system in closed-loop writes as: 
\begin{equation}
\label{pend-atay-gen}
    \Delta(\lambda; k_1,k_2,\tau_1,\tau_2):=\lambda^2-\frac{g}{l}+k_1e^{-\lambda\tau_1}+k_2e^{-\lambda\tau_2}.
\end{equation} 
Assume now that $(\tau_1,\tau_2)=(\tau,2\tau)$ with $\tau>0$, that is the commensurate delays case. Rescaling the time $t\mapsto t/\tau$, introducing the ``new'' parameter $\alpha^2=\tau^2 g/l\in \mathbb{R}_+$ and choosing $k_1=2\alpha^2$ and $k_2=-\alpha^2$, one gets:
\begin{equation}
\label{pend-atay-triple}
    \Delta(\lambda; \alpha):=\lambda^2-\alpha^2+2\alpha^2 e^{-\lambda}-\alpha^2 e^{-2\lambda}.
\end{equation}
It is easy to see that $\Delta$ in (\ref{pend-atay-triple}) can be factorized as $\Delta(\lambda;\alpha)=(\lambda-\alpha(1-e^{-\lambda}))(\lambda+\alpha(1-e^{-\lambda}))$. With no loss of generality, assume that $\alpha \in \mathbb{R}_+$. By taking into account the discussion done for the scalar system, we have the following observations:
\begin{itemize}
    \item First, for all $\alpha\in \mathbb{R}_+$, the characteristic function $\Delta$ in (\ref{pend-atay-triple}) has always an invariant root at the origin $\lambda=0$; its multiplicity is either $2$ or $3$. In fact, the multiplicity $3$ is reached if and only if $\alpha=1$. Furthermore, it is easy to show that the characteristic function $\Delta$ has three roots on the real axis.
    \item Second, if $\alpha=1$, excepting the triple root at the origin, all the remaining characteristic roots are all located in $\mathbb{C}_-$ and, thus, the root at the origin is dominant and the so-called \emph{multiplicity-induced-dominancy} still holds. 
    \item Finally, some basic but tedious algebraic manipulations allow concluding that for all $\alpha \in (0, 1)$, excepting the double root at the origin, all the remaining characteristic roots are all located in $\mathbb{C}_-$. Such a result suggests that the dominancy of the root at the origin is valid not only in the case when the maximal multiplicity is reached. However, such an observation is not generally true.
\end{itemize}

\begin{figure}[!ht]
\centering
\resizebox{0.9\columnwidth}{!}{\input{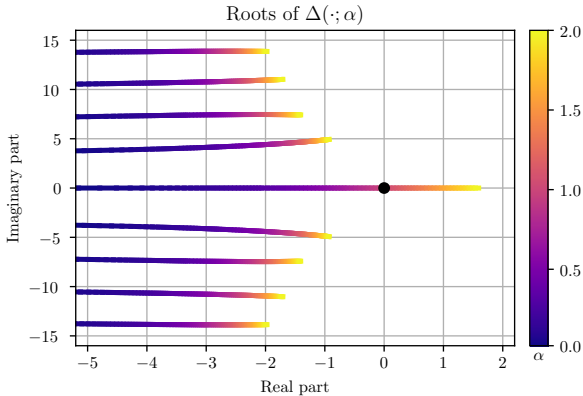}}
\caption{Roots of $\Delta(\cdot; \alpha)$ from \eqref{pend-atay-triple} for $\alpha \in [0, 2]$.}\label{fig:order-2}
\end{figure}

Similarly to Fig.~\ref{fig:order-1} for \eqref{eq:scalar-aux-char-f}, we represent in Fig.~\ref{fig:order-2} the behavior of the the roots of $\Delta(\cdot; \alpha)$ from \eqref{pend-atay-triple} as $\alpha$ varies in the interval $[0, 2]$, with the black dot representing the root at $\lambda = 0$. Due to the above factorization of $\Delta$, all roots of \eqref{eq:scalar-aux-char-f} for a given $\alpha$ are roots of both $\Delta(\cdot; \alpha)$ and $\Delta(\cdot; -\alpha)$, a fact that can be observed in Fig.~\ref{fig:order-2}. As before, the $\lambda = 0$ is the unique root in the case $\alpha = 0$ and, as $\alpha \to 0$, the real parts of all other roots of $\Delta(\cdot; \alpha)$ converge to $-\infty$. Moreover, still exploring the factorization of $\Delta$, one can check that, as in Remark~\ref{remk:root-locus}, any root of $\Delta(\cdot; \alpha)$ necessarily satisfies \eqref{eq:scalar-root-locus}, and, conversely, for every $\lambda \in \mathbb C$ satisfying \eqref{eq:scalar-root-locus} and not lying on the imaginary axis\footnote{Once again, the case of solutions of \eqref{eq:scalar-root-locus} on the imaginary axis can be retrieved in the limit $\alpha \to \pm \infty$.}, there exist exactly two values of $\alpha \in \mathbb R$ (a real number and its opposite) such that $\lambda$ is a root of $\Delta(\cdot; \alpha)$.

\begin{remark}
The planar inverted pendulum model (\ref{pend-atay}) was discussed by \cite{A1999}. More precisely, the author computed the stability regions of (\ref{pend-atay-gen}) in the parameter-space defined by the controller gains $(k_1,k_2)$ under the assumption that the delays are commensurate $\tau_1=\tau$ and $\tau_2=2\tau$. The idea to introduce delays in the control laws goes back to the 80s when \cite{SB1979} used the so-called ``proportional-minus-delay (PMD)'' controllers (see also \cite{SB1980}). The characterization of the codimension-three triple zero bifurcation of the inverted pendulum (\ref{InvertedPendulum}) by using various delay blocks including the ``PMD'' controller mentioned above can be found in \cite{BMN2015-SCL}. A deeper discussion on necessary conditions guaranteeing that multiple delay blocks may stabilize LTI SISO systems can be found in \cite{KNMM2005}. In the particular case of a chain of integrators, the explicit construction of the delay blocks can be found in \cite{NM2004}. For further discussions, we refer to \cite{mn14:siam} and the references therein. $\Box$ \end{remark}

\begin{remark}
Assume now that the planar inverted pendulum \eqref{pend-atay} is controlled by a standard delayed PD controller $u(t)=-k_p\theta(t-\tau)-k_d\dot \theta (t-\tau)$, then the delay margin $\tau_m$ guaranteeing the stability of closed-loop system for all $\tau \in [0, \tau_m)$ is $\tau_m=\sqrt{(2l)/g}$. This case study will be reconsidered in the forthcoming sections. For further discussions on such topics, we refer to \cite{Stepan2009} (see also \cite{Stepan89} and the references therein as well as \cite{A1999} for some further comparisons). $\Box$
\end{remark}

\section{Newton, Puiseux, Weierstrass and delay dynamics}
\label{sec:Pui}
Given an entire function $f(x,y)$, it is possible to reduce some of the analytic properties of $f$ to appropriate algebraic properties. To such a purpose, the following result (also known as the \emph{Weierstrass Preparation Theorem\/}) enables such a connection.
\begin{theorem}\label{theo:Weierstrass}
(\cite{Mailybaev:01}) Suppose that $f\left(z,\overrightarrow{p}\right)$ is an analytic function vanishing at the singular 
point $z_0\in\C$, $\overrightarrow{p}_0\in\C^{n}$, where $z=z_0$ is an $m$-multiple root of the equation $f\left(z,\overrightarrow{p}\right)=0$, i.e.,
\begin{eqnarray*}
f\left(z_0,\overrightarrow{p}_0\right)=:\left.\frac{\partial f}{\partial z}\right|_{(z_0,\overrightarrow{p}_0)}&=&\cdots=\left.\frac{\partial^{m-1}f}{\partial z^{m-1}}\right|_{(z_0,\overrightarrow{p}_0)}=0,\\ \left.\frac{\partial^{m}f}{\partial z^m}\right|_{(z_0,\overrightarrow{p}_0)}&\neq&0.
\end{eqnarray*}
\noindent
Then, there exists a neighborhood $U_0\subset\C^{n+1}$ of  $\left(z_0,\overrightarrow{p}_0\right)\in\C^{n+1}$ in which $f\left(z,\overrightarrow{p}\right)$ can be expressed as
\begin{equation}
\label{eq:Weierstrass}
f\left(z,\overrightarrow{p}\right)=W_p\left(z,\overrightarrow{p}\right)b\left(z,\overrightarrow{p}\right),
\end{equation}
\noindent
where $W_p\left(z,\overrightarrow{p}\right)=\left(z-z_0\right)^m+w_{m-1}\left(\overrightarrow{p}\right)\left(z-z_0\right)^{m-1}+\cdots+w_{0}\left(\overrightarrow{p}\right),$
and $w_{0}\!\left(\overrightarrow{p}\right)$,\ldots,$w_{m-1}\!\left(\overrightarrow{p}\right)$, $b\left(z,\overrightarrow{p}\right)$ are analytic functions uniquely defined 
by the function $f\left(z,\overrightarrow{p}\right)$ and $w_{i}\!\left(\overrightarrow{p}_0\right)=0$, $b\left(z_0,\overrightarrow{p}_0\right)\neq0$.
\end{theorem}
\begin{remark}
 The analytic function $W_p\left(z,\overrightarrow{p}\right)$ is known as the \textit{Weierstrass polynomial}. $\Box$
\end{remark}
\begin{remark}
It can be seen from Theorem \ref{theo:Weierstrass}, that since $b(z,\overrightarrow{p})$  is an holomorphic non-vanishing function at $(0,\overrightarrow{0})$,  then there must exists a neighborhood $\mathcal{O}(0,\overrightarrow{0})\subset\mathbb{C}^{n+1}$ at which $b(z,\overrightarrow{p})$ preserves the same property. Hence, 
based on this observation we can ensure that the root-locus of a given quasipolynomial $f=\Delta$ in the neighborhood $\mathcal{O}$ will be the same than 
the root-locus of $W_p(z,\overrightarrow{p})$. $\Box$
\end{remark}
\subsection{The Newton Diagram Method and Puiseux Series.}\label{sec:ndm}
\begin{figure}[!ht]
\centering
\includegraphics[clip=true,trim=36mm 5mm 92mm 9mm,scale=.34]{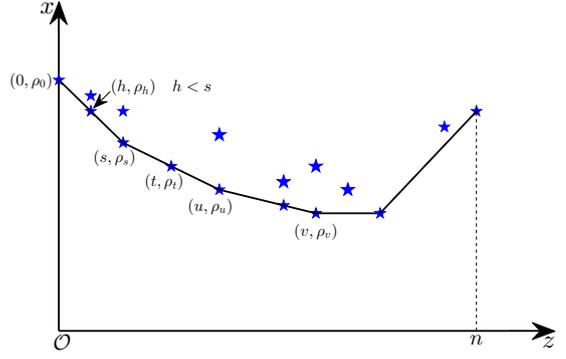}
\caption{Newton diagram for the pseudo-polynomial $f(x,y)$ given in \eqref{eqNewton}.}
\label{ddi}
\end{figure}

Given a known solution $\left(z_0,\overrightarrow{p}_0\right)$ of $f\left(z,\overrightarrow{p}\right)$, the local behaviour of the solution $z\left(\overrightarrow{p}\right)$ in 
the neighborhood $\C^{n}$ of $\overrightarrow{p}$ can be obtained by means of the \emph{Newton-diagram method}. 
Thus, in order to introduce such a procedure, let us consider the following notation (for more details, see, for instance, 
\cite{MM:19} and references therein). Let $f\left(x,y\right)$ be a \textit{pseudo-polynomial} in $y$, i.e.,
\begin{equation}
f\left(x,y\right)=\sum_{k=0}^{n}a_k(x)y^k,  \label{eqNewton}
\end{equation}
where the corresponding coefficients are given by:
\begin{equation}
a_k\left(x\right)=x^{\,\,\rho_k}\sum_{r=0}^{\infty}a_{rk}x^{r/q},\label{eqCNewton}
\end{equation}
and $a_{rk}\in \C$, $x$ and $y$ are complex variables,
$\rho_k\in \mathbb{Q}_+$, $q\in\mathbb{Z}_+$, $a_n(x)\not\equiv0$, and $a_0(x)\not\equiv0$. Then a solution of (\ref{eqNewton}) can be written in the form of a series as
\begin{equation*}
y=y_0+\alpha_1\left(x-x_0\right)^{\epsilon_1}+\alpha_2\left(x-x_0\right)^{\epsilon_2}+\cdots,
\end{equation*}
where $\epsilon_1$, $\epsilon_2$, $\ldots$, is an increasing
sequence of rational numbers. To determine the possible values of
$\epsilon_1$, $\alpha_1$, $\epsilon_2$, $\alpha_2$,
$\ldots$, it is necessary to consider the
\textit{Newton's diagram}. Since by simple translation, any point on a curve
can be moved to the origin, we will only consider expansions of the
solution of $f(x,y)=0$ around the origin. In this vein, we will
consider a solution of (\ref{eqNewton}) in the form:
\begin{equation}
y(x)=y_{\epsilon_1}
x^{\epsilon_1}+y_{\epsilon_2}x^{\epsilon_2}+y_{\epsilon_3}x^{\epsilon_3}+\cdots,
\label{solNewton}
\end{equation}
where $\epsilon_1<\epsilon_2<\epsilon_3<\cdots$,
$y_{\epsilon_1}\neq0$, or, in its compact form,
\begin{equation}
y(x)=y_{\epsilon_1}x^{\epsilon_1}+o \left(x^{\epsilon_1} \right). \label{VNewton}
\end{equation}
We have the following:
\begin{definition}[Newton's diagram and polygon]\label{def:ND}
Given a pseudo-polynomial of the form (\ref{eqNewton}) with
coefficients given by (\ref{eqCNewton}), plot $\rho_k$ versus $k$
for $k$ for $k=0,1,\ldots,n$ (if $a_k\left(\cdot\right)\equiv0$, the
corresponding point is disregarded). Denote each of these points by
$\pi_k=\left(k,\rho_k\right)$ and let
\[
\Pi=\left\{\pi_k\,:\,a_k(\cdot)\neq0\right\},
\]
be the set of all plotted points. Then, the set $\Pi$ will be called the Newton diagram, and the Newton polygon
associated with $f(x,y)$ will be given by the lower boundary of the convex hull of
the set $\Pi$.
\end{definition}
For a given pseudo-polynomial $f(x,y)$, Fig.~\ref{ddi}
simply illustrates Definition \ref{def:ND}.
The following result allows characterizing the solutions' structure of a given pseudo-polynomial (see, for instance, \cite{Wall:04}).
\begin{theorem}[Puiseux Theorem]
 \label{theo:Puiseuxtheo}
 The equation $f(x,y)=0$, with $f$ given in formal power series such that $f(0,0)=0$, possesses
 at least one solution in power series of the form:
 \[x=t^{q}, \quad y=\sum_{i=1}^{\infty}c_{i}t^{i}, \quad q\in\mathbb{N}.\]
\end{theorem}
\subsection{Asymptotic zero behavior characterization}
\label{sub-sec:NP}

The asymptotic behavior of the critical zeros of the quasipolynomial $f(\lambda,\tau)$ can be performed by means of the Newton diagram procedure.
To this end, since any critical solution $\left(\lambda^{\ast},\tau^{\ast}\right)$ can always be translated to the origin by appropriate shifts 
$\lambda\mapsto \lambda-\lambda^{\ast}$, $\tau\mapsto\tau-\tau^{\ast}$, hereinafter we will assume that $\left(\lambda^{\ast},\tau^{\ast}\right)=(0,0)$. Hence, 
for a $m$-multiple root $\lambda=0$ of $f$ at $\tau=0$, according to the Weierstrass Preparation theorem we will have that:
\begin{equation}
\label{eq:localf}
 f(\lambda,\tau)=\left(\lambda^m+w_{m-1}(\tau)\lambda^{m-1}\ldots+w_{0}(\tau)\right)b(\lambda,\tau).
\end{equation}
Now, with the aim of avoiding unnecessary computations, the following notation will be useful. For $i\in\llbracket 0,m-1\rrbracket$, denote by $n_i\in\mathbb{N}$ the first nonzero partial derivatives in $(\lambda,\tau)$ of $f$ at $(0,0)$, such that:
 \begin{subequations}\label{eq:ni}
\begin{eqnarray}
 \!\!\!\! f(0,0)=:\left.\frac{\partial^{i}f}{\partial \lambda^{i}}\right|_{(0,0)}&=&\cdots=\left.\frac{\partial^{i+n_{i}-1}f}{\partial \lambda^{i}\partial\tau^{n_{i}-1}}\right|_{(0,0)}=0,\\
       \left.\frac{\partial^{i+n_{i}}f}{\partial \lambda^{i}\partial\tau^{n_{i}}}\right|_{(0,0)}&\neq& 0.
\end{eqnarray}
\end{subequations}
The following result discloses the construction of the Newton diagram for a given quasipolynomial.
\begin{proposition}
\label{prop:ND}
 Let $\lambda=0$ be a $m$-multiple root at $\tau=0$ of the quasipolynomial $\Delta(\lambda,\tau)$, and assume that $n_{0}<\infty$. Then, the Newton 
 diagram of $f$ at $(0,0)$ is given by $\Pi=\left\{(0,n_{0}),\dots,(m-1,n_{m-1}),(m,0)\right\}$.
\end{proposition}
\begin{example}\label{exmp:1a}
 To illustrate the previous result, consider the quasipolynomial (borrowed from \cite{Cai:12}):
\begin{equation}
\label{eq:exp2:sec4}
 \Delta(\lambda,\tau)=-P_0\left(\lambda\right)+P_1\left(\lambda\right)e^{-s\tau}+e^{-2s\tau},
\end{equation}
where:
\begin{equation*}
 P_0\!\left(\lambda\right):=\frac{\pi}{2}\lambda^5+\frac{\pi}{2}\lambda^3+\lambda^2,\,
 P_1\!\left(\lambda\right):=\frac{\pi}{2}\lambda^3-\lambda^2+\frac{\pi}{2}\lambda+1,
\end{equation*}
and with $\lambda=j$ a root at $\tau=\pi$ of multiplicity $m=3$. First, we derive the constants $n_i$ considered in (\ref{eq:ni}):
\begin{align*}
 \!\!\!\!\!\left.\frac{\partial \Delta}{\partial\tau}\right|_{(j,\pi)}\!\!\!\!\!&=0,\,\left.\frac{\partial^{2}\Delta}{\partial\tau^2}\right|_{(j,\pi)}=-2\,\,\,\Rightarrow n_0=2,\\ \quad \left.\frac{\partial^{2}\Delta}{\partial \lambda\partial\tau}\right|_{(j,\pi)}
             \!\!\!\!\!&=2+j\pi,\Rightarrow n_1=1,\\
 \left.\frac{\partial^{3}\Delta}{\partial \lambda^{2}\partial\tau}\right|_{(j,\pi)}\!\!\!\!\!&=-\left(5\pi+j(4\pi^{2}+6)\right)\Rightarrow n_2=1,\\
        \left. \frac{\partial^{3}\Delta}{\partial \lambda^{3}}\right|_{(j,\pi)}\!\!\!\!\!&=-3\pi(-6-5j\pi+\pi^{2}).
\end{align*}
 Summarizing, we have $\left(n_0,n_1,n_2\right)=(2,1,1)$. According to Proposition \ref{prop:ND}, we have: $\Pi=\left\{(0,2),\,(1,1),(2,1),(3,0)\right\}$. 
 Such points are depicted in Fig.~\ref{Fig:1b}.
 \begin{figure}[!ht]
\centering
\includegraphics[clip=true,trim=100mm 15mm 70mm 2mm,scale=0.325]{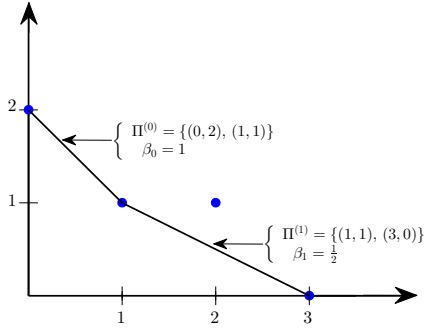}
\caption{Newton Diagram for the quasipolynomial \eqref{eq:exp2:sec4}.}\label{Fig:1b}
\end{figure}
\end{example}

\begin{remark}
Note that it is possible to have some $\kappa\in\mathbb{N}$ for which $n_{0}=n_{1}=\cdots=n_{\kappa-1}=\infty$. Then, under this situation, the Newton diagram method can not be applied directly. However, it is worth noting that since $w_i$ are analytic functions, the previous situation is equivalent to $w_{i}(\tau)\equiv 0$ for $0\leq i\leq \kappa-1$. Hence, $f$ will be locally given by

\[
\!f(\lambda,\tau)\!=\!\lambda^{\kappa}\!\!\left[\lambda^{m-\kappa}\!\!+w_{m-\kappa}(\tau)\lambda^{m-\kappa-1}\!\!+\cdots+w_{\kappa}(\tau)\right]\!b(\lambda,\tau).
\]
Thus, there are $\kappa$-\textit{invariant} solutions at $\lambda=0$ for all $\tau$ and $m-\kappa$ solutions of the form
\[
\lambda_{i}(\tau)=\sum_{\ell=1}^{\infty}c_{i,\ell}\,\tau^{\ell/m_{i}},
\]
where $m_i<m$. Moreover, under this consideration the Newton polygon will be given by $\Pi=\left\{(\kappa,n_{\kappa}),\ldots,(m,0)\right\}$.
If such number $\kappa$ does not exist (i.e., if such situation does not happen), then $\kappa$ will be simply defined as $\kappa:=0$. $\Box$
\end{remark}

\subsection{Puiseux algorithm \& Weierstrass polynomial}
\label{sec:Weier}

As seen in Section~\ref{sec:Pui} from Theorem \ref{theo:Weierstrass}, we know that $b(z,\overrightarrow{p})$  is a holomorphic non-vanishing function at $(0,\overrightarrow{0})$. Then, there must exist some neighborhood $\Omega(0,\overrightarrow{0})\subset\C^{n+1}$ at which $b(z,\overrightarrow{p})$ preserves the same property. Hence, based on this observation we can ensure that the root-locus of a given quasipolynomial $\Delta$ in the neighborhood $\Omega$ will be the same as the root-locus of $W_p(z,\overrightarrow{p})$. Thus, bearing this observation in mind, the next result allows computing the Weierstrass polynomial.

\begin{proposition}
\label{prop:Weierstrass}
 Let $n_{i}<\infty$ be defined as in \eqref{eq:ni} such that 
$n_{0}>n_{1}>\cdots>n_{m-1}$. Then, the coefficients $w_{i}(\tau)$ of the associated Weierstrass polynomial $W$ have order $\ord{w_{i}(\tau)}{}=n_{i}$. Moreover, the leading 
 terms are given by:
 \[
  w_{i}(\tau)=\left(\frac{m!}{i!n_{i}!\left.\frac{\partial^mf}{\partial \lambda^m}\right|_{(0,0)}}\left.\frac{\partial^{i+n_{i}}f}{\partial\tau^{n_{i}}\partial \lambda^{i}}\right|_{(0,0)}\right)\tau^{n_{i}}
              +o(\tau^{n_{i}}),
 \]
 for $i \in \llbracket 0,m-1\rrbracket$.
\end{proposition}
As discussed in Section~\ref{sec:Pui}, Theorem \ref{theo:Puiseuxtheo} allows revealing the solution structure of a given quasipolynomial $\Delta$. In fact, the equation $\Delta(\lambda;\tau)=0$ defines a manifold $\mathcal{C}\in \mathbb{C}^{2}$ which is composed by the finite union of $r$ branches $\lambda_{i}\left(\tau^{1/m_i}\right)$, each of these branches can be expressed as a Puiseux series:
\begin{equation}
\label{def:Puiseux}
 \lambda_{i\sigma_i}(\tau)=c_{i\sigma_i}\tau^{\frac{1}{m_i}}+o\left(|\tau|^{\frac{1}{m_i}}\right),
\end{equation}
for $i \in \llbracket 0,r-1 \rrbracket$, $\sigma_i \in \llbracket 1,m_i\rrbracket$, where each branch has multiplicity $m_i$, such that $m=m_1+m_2+\cdots+m_{r}$. In the case when $r=1$, then $\lambda_{i\sigma_i}$ and $c_{i\sigma_i}$  will be simply denoted by $\lambda_{\sigma_i}$ and $c_{\sigma_i}$, respectively.\\ 
Now, with the aim to classify the nature of any solution, let us introduce the following definition.
\begin{definition}
 We say that there is a \emph{Complete Regular Splitting} (CRS) property of the solution $s^{\ast}=0$ at $\tau^{\ast}=0$ if 
 $c_{i\sigma_i}\neq0,\,\,\forall i$. For the \emph{Regular Splitting} (RS) property, some of the coefficients $c_{i\sigma_i}$ for which 
 $m_i=1$ may be equal to zero. In the remaining cases of the coefficient $c_{i\sigma_i}$ we say that \textit{Non-Regular Splitting} (NRS) property is present.
\end{definition}
To deal with the splitting properties of a given solution, we will consider the Newton diagram method in conjunction with the Weierstrass polynomial and the Puiseux Theorem. To such an end, Algorithm \ref{alg} will be extremely useful.
\begin{algorithm}[h!]
\caption{Puiseux Series Expansion Algorithm.}
\label{alg}
Let $f(\lambda;\tau)$ have a critical pair such that $\lambda^{\ast}=j\omega^{\ast}$ is a $m$-multiple root at $\tau=\tau^{\ast}$. Consider the initial 
values as $r:=0$, $i_{-1}:=\kappa$ and $\ell_{-1}:=n_\kappa$.
\DontPrintSemicolon

\While{$i_{r-1}<m$}{
set $\mathcal{E}_{r}:=\left\lbrace \frac{\ell-\ell_{r-1}}{i_{r-1}-i}:\left(i,\ell\right)\in\Pi,\text{ and }i>i_{r-1}\right\rbrace$;\;

let $\beta_{r}:=\max\mathcal{E}_{r}$ and \;

$\Pi^{(r)}:=\left\{\left(i,\ell\right)\in\Pi:\beta_r\equiv\frac{\ell-\ell_{r-1}}{i_{r-1}-i}\right\}\cup\left\{\left(i_{r-1},\ell_{r-1}\right)\right\}$;\;

set $(i_{r},\ell_{r})\in\Pi^{(r)}$ such that $i_r\geq i,\,\,\forall(i,\ell)\in\Pi^{(r)}$;\;

set $m_r:=i_r-i_{r-1}$ and $r=r+1$.\;

}
\end{algorithm}

The following result allows splitting identification:
\begin{proposition}\label{prop:Splitting}
 Let $\lambda^{\ast}=j\omega^{\ast}$ at $\tau=\tau^{\ast}$ be a $m$-multiple critical root of the quasipolynomial $\Delta(\lambda;\tau)$. Assume that 
$r$, $\beta_i$, $\left(i_i,\ell_i\right)$, $m_i$ and $\Pi^{(i)}$, for $i \in \llbracket 0,r-1\rrbracket$ are given by the Algorithm \ref{alg}. Then the following 
 properties hold:
\begin{description}
  \item[(i)] if $m_i\cdot\beta_{i}\equiv1$, $\forall i\in \llbracket 0,r-1\rrbracket$, then the solution $\left(j\omega^{\ast},\tau^{\ast}\right)$
  of $\Delta(\lambda;\tau)$ has the CRS property;
  \item[(ii)] if some $\beta_{i}$ satisfies $m_i\cdot\beta_{i}>1$ for $m_{i}>1$, then NRS  property for the solution $\left(j\omega^{\ast},\tau^{\ast}\right)$ occurs;
  \item[(iii)] if the pairs $\left(m_k,\beta_k\right)$ that do not fulfill (i), satisfy the inequality 
  $\beta_k\geq m_k\equiv1$, then the solution $\left(j\omega^{\ast},\tau^{\ast}\right)$ of $\Delta(\lambda,\tau)$ has the RS property.
\end{description}

\begin{corollary}
\label{cor:n0}
With the hypothesis above (Proposition \ref{prop:Splitting}), assume that $n_0=1$. Then at $\tau=\tau^{\ast}$, the $m$-roots of $\Delta(\cdot;\tau)$ have the CRS property and can be expanded as:
 \begin{equation}
  \lambda_{\sigma_i}\left(\tau\right)=j\omega^{\ast}+c_{\sigma_i}\left(\tau-\tau^{\ast}\right)^{\frac{1}{m}}+
                            o\left(\left|\tau-\tau^{\ast}\right|^{\frac{1}{m}}\right),
 \end{equation}
 for $\sigma_i \in \llbracket 1,m\rrbracket$. Moreover, the following properties hold:
 \begin{description}
  \item[(i)] if $m=2$ and $\Re\left(c_{sigma_i}\right)\neq0$ with $\sigma_i\in\{1,2\}$. Then for $\tau>\tau^{\ast}$ sufficiently close 
             to $\tau^{\ast}$, one of the zeros $\lambda_{\sigma_i}\left(\tau\right)$ will enter $\mathbb{C}_+$, whereas the other one will enter $\mathbb{C}_-$;
  \item[(ii)] if $m>2$, then at least one of the zeros $\lambda_{\sigma}\left(\tau\right)$ will enter $\mathbb{C}_+$.
 \end{description}
\end{corollary}
\end{proposition}
As mentioned earlier, the Weierstrass polynomial will allow us to analyze the stability behavior of the imaginary characteristic roots. In this vein, we have the following:
\begin{proposition}\label{prop:cdr}
 Let $n_0<\infty$ and $\lambda^{\ast}=j\omega^{\ast}$ be a $m$\nobreakdash-mul\-ti\-ple root of $\Delta(\lambda;\tau)$ at $\tau=\tau^{\ast}$. Assume that 
 $r$, $\beta_i$, $\left(i_i,\ell_i\right)$, $m_i$ and $\Pi^{(i)}$, for $i \in \llbracket 0,r-1\rrbracket$ are given by the Algorithm \ref{alg}. 
 Then, at  $\tau=\tau^{\ast}$ the $m$-zeros of $\Delta$ can be expanded as 
 \begin{equation}
   \lambda_{i\sigma_i}(\tau)=j\omega^{\ast}+c_{i\sigma_i}\left(\tau-\tau^{\ast}\right)^{\beta_{i}}+
                       o\left(\left|\tau-\tau^{\ast}\right|^{\beta_{i}}\!\right),\label{eq:sj}
 \end{equation}
 for $i \in \llbracket 0,r-1\rrbracket$, $\sigma_i \in \llbracket 1,m_{i}\rrbracket$ and $m=m_0+\cdots+m_{r-1}$. 
 Where $c_{i\sigma_i}$ are roots of the polynomial $\mathcal{P}_{i}:\mathbb{C}\mapsto \mathbb{C}$,
 \begin{equation}
 \label{eq:pnu}
  \mathcal{P}_{i}(z):=\sum_{k=i_{i-1}}^{i_{i}}\,a_{k,0}z^{k-i_{i-1}},\quad \text{s.t. }\left(k,\eta_{k}\right)\in\Pi^{(i)},
 \end{equation}
 where the coefficient $a_{k,0}\in\C$ is given by
 \begin{equation}
 \label{eq:acoef}
  a_{k,0}=\left(\frac{m!}{k!\eta_{k}!\left.\frac{\partial^m\Delta}{\partial \lambda^m}\right|_{(0,0)}}\right)\left.\frac{\partial^{k+\eta_{k}}\Delta}{\partial\tau^{\eta_{k}}\partial \lambda^{k}}\right|_{(0,0)}.
 \end{equation}
Furthermore, for  $\tau>\tau^{\ast}$ sufficiently close to $\tau^{\ast}$, the zeros $\lambda_{i\sigma_i}(\tau)$ will enter $\mathbb{C}_+$ (or $\mathbb{C}_-$) if
 \begin{equation}
  \Re\{c_{i\sigma_i}\}>0 \;(<0).\label{CDR}
 \end{equation}
\end{proposition}

The following results allow a further characterization:

\begin{proposition}
\label{prop:taylor}
Let $\lambda^{\ast}=j\omega^{\ast}$ be a $m$-multiple root of $\Delta(\lambda;\tau)$ at $\tau=\tau^{\ast}$. Assume that 
$\beta_i$ and $m_i$ for $i \in \llbracket \kappa,r-1\rrbracket$ are given by the Algorithm \ref{alg}. If $\beta_{i}=1$, then the following statements hold:
 \begin{itemize}
  \item[(i)] the equation $\Delta(\lambda;\tau)=0$ has $m_{i}$-solutions of the form
           \begin{equation}
            \lambda_{i\sigma_i}(\tau)=j\omega^{\ast}+c_{i\sigma_i}\left(\tau-\tau^{\ast}\right)+o\left(\left|\tau-\tau^{\ast}\right|\right), \label{eq:taylor1}
           \end{equation}
 	   with $\sigma_i\in\llbracket 1,m_{i}\rrbracket$ and where $c_{i\sigma_i}$ is a root of the polynomial $\mathcal{P}_i$ defined in \eqref{eq:pnu};
 \item[(ii)] if $c_{i\sigma_i}$ is a simple root of $\mathcal{P}_i$ then, there are $m_{j}$-solutions
           expanded as a Taylor series in the form
           \[
            \lambda_{i\sigma_i}(\tau)=j\omega^{\ast}+c_{i\sigma_i}\left(\tau-\tau^{\ast}\right)+c_{i\sigma_i}^{(1)}\left(\tau-\tau^{\ast}\!\right)^{1+\beta_i^{(1)}}+\cdots,
           \]
           where $\beta_i^{(1)}\in\mathbb{N}$.
 \end{itemize}
\end{proposition}

\begin{proposition}
\label{prop:puiseux}
 Let $\lambda^{\ast}=j\omega^{\ast}$ be a $m$-multiple root of $\Delta(\lambda,\tau)$ at $\tau=\tau^{\ast}$. Assume that 
 $\beta_h$, $m_h$ and $\left(i_{h},\ell_{h}\right)\in\Pi^{(h)}$ for $h \in \llbracket 0,r-\kappa-1\rrbracket$ are given by the Algorithm \ref{alg}. If $\beta_{h}=1/m_{h}$, then 
 $\Delta(\lambda;\tau)=0$ has $m_{u}$-solutions given by
 \[
  \lambda_{h\sigma_h}(\tau)=j\omega^{\ast}+c_{h}\Theta_{\sigma_h}\left(\tau-\tau^{\ast}\right)^{1/m_{h}}+o(\left|\tau-\tau^{\ast}\right|^{1/m_{h}}),
 \]
 with $\sigma_h\in\llbracket1,m_{j}\rrbracket$ where $\Theta_{\sigma_h}=\exp\left(j\frac{\theta_h+2\pi (\sigma_h-1)}{m_{h}}\right)$, $\theta_h=\arg(c_{h}^{m_{h}})$ and 
 $c_{h}=\left|a_{i_{h-1},0}/a_{i_{h},0}\right|^{1/m_{h}}$.
\end{proposition}

\begin{example}[Inverted pendulum]
Reconsider the stabilization of the planar inverted pendulum without friction (\ref{pend-atay-gen}) under the assumption of commensurate delays $(\tau_1,\tau_2)=(\tau,2\tau)$, with $\tau \in \mathbb{R}_+$. Thus, the characteristic function rewrites as:
\begin{equation}
\label{eq:atay}
\Delta(\lambda;\tau)= \lambda^2 -\frac{g}{l}+k_1e^{-\lambda\tau}+ k_2e^{-2\lambda\tau}.
\end{equation}
By setting $k_1+k_2=g/l$, we have that $\Delta(\lambda,0)=\lambda^2$ and $\Delta(0;\tau)=0$, for all $\tau\in \mathbb{R}_+$. Moreover, if $k_1=-2k_2$ the first partial derivative $\partial_\lambda\Delta$ evaluated at $\lambda=0$ also vanishes.\\ 
Consider the delay interval $0<\tau<\sqrt{2 \abs{k_1 + 4 k_2}}$, the root at the origin $\lambda^\ast=0$ has multiplicity $m=2$. The corresponding Weierstrass polynomial writes as $W_p(\lambda;\tau)=\lambda^2+w_1(\tau)\lambda+w_0(\tau)$. For its computation, the first partial derivatives of $\Delta$ at $\lambda^\ast=0$ and $\tau^\ast=\sqrt{2(k_1+4k_2)}$ are:
  \begin{equation}\label{eq:atay_n_a}
\left\{  \begin{array}{l} \left.\,\,\,\,\,\,\,\frac{\partial^{n_0} \Delta}{\partial\tau^{n_0}}\right|_{\lambda^\ast=0}=0, \forall n_0\in\mathbb{N} \Rightarrow n_0=\infty,\vspace{5mm} \\
 \left.\frac{\partial^{n_1+1} \Delta}{\partial\tau^{n_1}\partial \lambda}\right|_{\lambda^\ast=0}=0, \forall n_1\in\mathbb{N} \Rightarrow n_1=\infty,
 \end{array}\right.
  \end{equation}
  implying that $\kappa=2$. In the light of \eqref{eq:atay_n_a}, the coefficients $w_0=w_1\equiv0,\,\forall\tau\in\mathbb{R}_+$. Thus, there are two-invariant solutions at $\lambda=0$ and, around the origin, $\Delta$ can be written as $\Delta(\lambda,\tau)\equiv\lambda^2\widehat{\Delta}(\lambda,\tau)$.\\
  With the parameters choice $k_2=-g/l$, $k_1=2 g/l$ and $\tau^\ast=\sqrt{l/g}$, the multiplicity of $\lambda^\ast=0$ is $m=3$. Thus, following the Weierstrass Preparation Theorem, the corresponding local behavior is captured by
$\Delta(\lambda;\tau)=\left[\lambda^3+w_2(\tau)\lambda^2+w_1(\tau)\lambda+w_0(\tau)\right]b(\lambda,\tau)$.    
  Similarly to the double root case \eqref{eq:atay_n_a}, $n_0=n_1=\infty$, 
  leaving the following derivatives to be determined:
  \begin{equation}\label{eq:ata_n_b}
 \left\{  \begin{array}{l}\left.\frac{\partial^3 \Delta}{\partial\tau\partial \lambda^2}\right|_{(\lambda^\ast,\tau^\ast)}=-4\left(\sqrt{g/l}\right)^3 \Rightarrow n_2=1,\\
 \left.\frac{\partial^3 \Delta}{\partial \lambda^3}\right|_{(\lambda^\ast,\tau^\ast)}= 6\sqrt{l/g}.
 \end{array}\right.
  \end{equation}
  Therefore, the Weierstrass coefficient $w_2\neq 0$, meaning that $\kappa=2$. In other words, we have 2-invariant solutions at $\lambda=0$ and a solution $\lambda(\tau)=-w_2(\tau)$. The first term in the expansion of $w_2$ can be obtained by using \eqref{eq:ata_n_b} resulting in the following Weierstrass polynomial
$W_p(\lambda;\tau)=\lambda^2\left(\lambda+w_2(\tau)\right)$.
  Therefore, two roots $\lambda_{1,2}\equiv0$ remain invariant under delay variations, and the third root given by $\lambda_3(\tau)=-w_2(\tau)$. With this results, the Newton polygon has a single segment with two points as shown in Table \ref{tab:atay}.
\begin{table}[h!]
 \centering
 \caption{Summary: Inverted pendulum \eqref{eq:atay}.}
  {\begin{tabular}{l l l} 
  \toprule
   Initial Data & Algorithm Output & $\{z\in\mathbb{C}:\mathcal{P}_j(z)=0\}$ \\ 
  \midrule
   $m=3$, $n_{2}=1$ & $r=1, m_{0}=1, \beta_{0}=1$ & $\mathcal{P}_{0}(z):=z+2g/l$ \\
   $\Pi\!=\!\left\{(2,1),(3,0)\right\}$ & $\Pi^{^{(0)}}=\{(2,1),(3,0)\}$ & $\{c_{0,1}=-2g/l\}$ \\ 
  \bottomrule
  \end{tabular}}\label{tab:atay}
 \end{table}
  Now, the asymptotic behavior of the of the solution $\lambda_3$ is obtain by means of Proposition \ref{prop:cdr} Table \ref{tab:atay}, resulting in
   $\lambda_3(\tau)=-2g\frac{g}{l}\tau+o(\tau)$.
  as depicted in Fig.~\ref{fig:atay}.
  \begin{figure}[!ht]
 \centering  
  \includegraphics[trim=28mm 5mm 25mm  10mm,clip,width=\columnwidth]{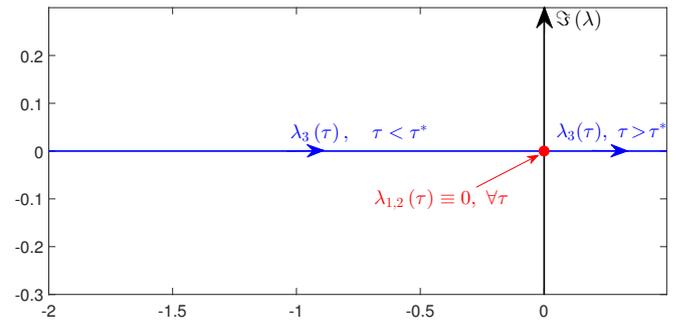}
 \caption{Root locus behavior ($2$-invariant root at $\lambda^\ast=0$) for \eqref{eq:atay} with $l=2g$ ($k_1=1$, $k_2=-0.5$)}\label{fig:atay}
 \end{figure}
\end{example}

\section{Frequency-sweeping curves and imaginary characteristic roots}
\label{sec:fsc}

Consider now the commensurate delay case and we are interested on the effects induced by the delay parameter. The corresponding characteristic function writes as:
\begin{equation}
\label{generaltimedelaysystemofRetardedtype}
\Delta (\lambda; \tau) = P_0(\lambda)+\sum_{i}^{n_d}P_i(\lambda)e^{-i \lambda \tau},   
\end{equation}
where $P_i\in \mathbb{R}[x]$, for all $i \in \llbracket0, n_d\rrbracket$ with $\deg(P_0)>\max_i(\deg(P_i))$, for all $i \in \llbracket1, n_d\rrbracket$. In this frame, for a critical pair $({\lambda _\alpha },{\tau _{\alpha ,k}})$,
denote by $n\in \mathbb{N}_+$ the multiplicity of $\lambda _\alpha$
at $\tau _{\alpha ,k}$. Clearly, a critical imaginary root is called
a \emph{simple critical imaginary root} (a \emph{multiple critical imaginary
root}) if the
corresponding index $n=1$ ($n>1$). In other words, the index
$n$ simply implies that for $\lambda=\lambda
_\alpha$ and $\tau=\tau _{\alpha ,k}$,
\begin{eqnarray}
\label{definitionofn} {\Delta_{{\lambda^0}}} =  \cdots  = {\Delta_{{\lambda^{n - 1}}}} = 0,\, {\Delta_{{\lambda^n}}} \ne 0.
\end{eqnarray}
Next, introduce the index $g\in \mathbb{N}_+$ at
$(\lambda _\alpha,\tau _{\alpha ,k})$, by which we may
\emph{artificially} treat $\tau _{\alpha ,k}$ as a $g$-multiple root
for $\Delta (\lambda; \tau)=0$ when $\lambda=\lambda _\alpha$, having the
property that when $\lambda=\lambda _\alpha$ and $\tau=\tau _{\alpha
,k}$,
\begin{eqnarray}
\label{definitionofg} {\Delta_{{\tau ^0}}} =  \cdots  = {\Delta_{{\tau ^{g -
1}}}} = 0,\,{\Delta_{{\tau ^g}}} \ne 0.
\end{eqnarray}
Suppose that $(\alpha, \beta)$ (with $\beta>0$) is a critical pair with the index $n$. Near this critical
pair, there exist $n$ roots $ {\lambda_i}(\tau )$  continuous w.r.t.\ $\tau$ satisfying $\alpha = {\lambda_i}(\beta )$, $i = 1,
\ldots, n$. Under some perturbation $\varepsilon
$ ($ - \varepsilon $) on $\beta$, the $n$ roots are expressed by $ {\lambda_i}(\beta +
\varepsilon)$ ($ {\lambda_i}(\beta - \varepsilon)$), $i = 1, \ldots
,n$. Denote the number of unstable roots among ${\lambda_1}(\beta +
\varepsilon ), \ldots ,{\lambda_n}(\beta  + \varepsilon )$
(${\lambda_1}(\beta - \varepsilon ), \ldots ,{\lambda_n}(\beta  -
\varepsilon )$) by $N{U_\alpha }({\beta ^ + })$ ($N{U_\alpha
}({\beta ^ - })$). With these notations, we define:
\begin{equation}\label{DeltaNU}
\mathrm{Var}(N{U_\alpha }(\beta )) \buildrel \Delta \over =
N{U_\alpha }({\beta ^ + }) - N{U_\alpha }({\beta ^ - }).
\end{equation}
The notation $\mathrm{Var}(N{U_\alpha }(\beta ))$  stands for the change of $NU(\tau)$ caused by the variation of the
critical imaginary root $ \lambda = \alpha  $ as $\tau$ increases from $\beta - \varepsilon$ to $\beta
+ \varepsilon$.
\begin{remark}
To the best of the authors' knowledge, the first systematic discussion on the number of unstable roots by using the continuity of the roots with respect to the delay parameter can be found in \cite{kashiwagi:65}, where the author introduced a similar concept to ``$(NU)$'' called \emph{stability indicative function}. $\Box$
\end{remark}

We here give the procedure for generating the frequency-sweeping curves. The characteristic function $\Delta (\lambda; \tau)$ can be transformed by letting $ z = e^{ - \tau \lambda} $ into a two-variate (auxiliary) polynomial $P_a$:
\begin{equation}
\label{lambdaz} P_a(\lambda,z) =\sum\limits_{i = 0}^{n_d} {P_i
(\lambda)z^i }.
\end{equation}

\emph{Frequency-Sweeping Curves (FSCs)}: Sweep $\omega  \ge 0$ and
for each $\lambda=j \omega$ we have $q$ solutions of $z$ such that
$P_a (j \omega,z)=0$ (denoted by ${z_i}(j\omega ), i \in \llbracket 1,n_d \rrbracket$) leading to $n_d$  frequency-sweeping
curves $\Gamma_i(\omega)$: $\left| {{z_i}(j\omega )} \right|$ vs.
$\omega$, $i \in \llbracket 1,n_d \rrbracket$.

If $(\lambda_\alpha, \tau_{\alpha,k})$ is a critical pair with index
$g$, ``$g$'' FSCs intersect $ \Im _1 $\footnote{We denote  by $ \Im _1 $ the line parallel to the abscissa axis with ordinate equal to one.} at $\omega=\omega_\alpha$ and the frequency $\omega_\alpha$ is called a \emph{critical
frequency}.

In the sequel, we introduce some necessary notations concerning the asymptotic behavior of FSCs. For a deeper discussion, we refer to Chapter 8 of \cite{lnc:15-briefs}. Suppose $({\lambda _\alpha },{\tau _{\alpha ,k}}), k \in
\mathbb{N}$, is a set of critical pairs (as usually assumed,
${\lambda _\alpha } \ne 0$) with the index $g$ ($g$ is a constant w.r.t.\ different $k$, see Property 1.2 of \cite{lnc:15-briefs}). There must exist ``$g$'' FSCs such that ${z_i}(j{\omega _\alpha}) = {z_\alpha}= {e^{ -
{\tau _{\alpha ,0}}{\lambda _\alpha }}}$ intersecting ${\Im _1}$
when $\omega  = {\omega _\alpha }$. Among such $g$
FSCs, we denote the number of the
FSCs when $\omega = {\omega _\alpha} +
\varepsilon $ ($\omega = {\omega _\alpha} - \varepsilon $) above the line $
\Im _1 $ by $N{F_{{z_\alpha}}}({\omega _\alpha} + \varepsilon )$
($N{F_{{z_\alpha}}}({\omega _\alpha} - \varepsilon )$). Introduce
now a new notation $\mathrm{Var}( N{F_{{z_\alpha}}})({\omega _\alpha})$ as
\begin{eqnarray}
\label{definitionDeltaNF} \mathrm{Var}(N{F_{{z_\alpha}}})({\omega _\alpha})
\buildrel \Delta \over = N{F_{{z_\alpha}}}({\omega _\alpha} +
\varepsilon ) - N{F_{{z_\alpha}}}({\omega _\alpha} - \varepsilon ).
\end{eqnarray}

\begin{theorem}[\cite{lnc:17-ieee}]
\label{InvarianceTheorem} Let ${\lambda
_\alpha }$ be a critical imaginary root of the characteristic function
(\ref{generaltimedelaysystemofRetardedtype}). Then 
$\mathrm{Var}(N{U_{{\lambda _\alpha}}})({\tau _{\alpha,k}})$ is a constant
$\mathrm{Var}(N{F_{{z_{\alpha}}}})({\omega _\alpha})$ for all ${\tau
_{\alpha ,k}}
> 0$.
\end{theorem}

The contribution of Theorem \ref{InvarianceTheorem} is twofold:
\begin{itemize}
\item First, it provides a simple method (observing the FCS) to compute $\mathrm{Var}(N{U_{{\lambda _\alpha }}})({\tau _{\alpha ,k}})$, without invoking the Puiseux series.
\item Second, a very interesting \emph{invariance property} is claimed: \emph{For a critical imaginary root ${\lambda
_\alpha }$, the asymptotic behavior has the same effect on the stability (more precisely, on $NU(\tau)$) at all the corresponding positive critical delays ${\tau
_{\alpha ,k}} > 0$}.
\end{itemize}
With the \emph{invariance property}, we can overcome the peculiarity that a critical imaginary root has infinitely many critical delays.

\begin{remark}
By using a different argument, \cite{jm:10} addressed the invariance property in the case when $n=2$ and $g=1$. $\Box$
\end{remark}

\begin{example}\label{MotivatingExampleVariantMultiplicityTACFP2017} Consider
a time-delay system with the characteristic function $\Delta (\lambda; \tau) = \sum\nolimits_{i = 0}^4
{{P_i}(\lambda ){e^{ - i\tau \lambda }}}$ where ${P_0}(\lambda ) =
\frac{{15}}{8}{\pi ^2}{\lambda ^6} + (\frac{{11}}{4}\pi  -
\frac{{15}}{8}{\pi ^2}){\lambda ^4} + \frac{9}{2}\pi {\lambda ^3} +
(1 + \frac{1}{2}\pi  - \frac{{75}}{8}{\pi ^2}){\lambda ^2} + (3 +
\frac{9}{2}\pi )\lambda + 1 - \frac{9}{4}\pi  - \frac{{45}}{8}{\pi
^2}$, ${P_1}(\lambda ) = \frac{5}{4}\pi {\lambda ^5} +
\frac{{11}}{2}\pi {\lambda ^4} + (1 + \frac{7}{2}\pi ){\lambda ^3} +
(\pi  + 7){\lambda ^2} + (11 + \frac{9}{4}\pi )\lambda  + 4 -
\frac{9}{2}\pi $,  ${P_2}(\lambda ) = \frac{5}{4}\pi {\lambda ^5} +
\frac{{11}}{4}\pi {\lambda ^4} + (3 - \pi ){\lambda ^3} + (13 +
\frac{1}{2}\pi ){\lambda ^2} + (15 - \frac{9}{4}\pi )\lambda  + 6 -
\frac{9}{4}\pi $, ${P_3}(\lambda ) = 3{\lambda ^3} + 9{\lambda ^2} +
9\lambda  + 4$, and ${P_4}(\lambda ) = {\lambda ^3} + 2{\lambda ^2}
+ 2\lambda  + 1$.\\
We study the asymptotic behavior of critical pairs $(j,(2k + 1)\pi )$, with $g=2$.
The frequency-sweeping curves are given in Fig.~\ref{FSCExampleInvarianceConfirmationofMotivatingExampleVariantMultiplicityTACFP2017}. According to Theorem \ref{InvarianceTheorem}, we know from Fig.~\ref{FSCExampleInvarianceConfirmationofMotivatingExampleVariantMultiplicityTACFP2017} that $\Delta N{U_j}((2k + 1)\pi ) = 0$ for all $k \in \mathbb{N}$.\\
In fact, the  asymptotic behavior of critical pairs $(j,(2k + 1)\pi )$ is complicated. The
multiplicity $n$ of the critical imaginary root $\lambda = j$ is 2, 3, 4, 2, when $\tau$ is $\pi$, $3\pi$, $5 \pi$, $7 \pi$, respectively. The  Puiseux series are
all degenerate:
\begin{eqnarray*} 
&& \left\{
{\begin{array}{*{20}{l}}
\begin{array}{l}
\delta \lambda  = 0.1592j\delta \tau  + (0.5371 - 0.3138j){(\delta
\tau )^2}\\
\quad\quad\quad  + o({(\delta \tau )^2}),
\end{array}\\
\begin{array}{l}
\delta \lambda  = 0.0796j\delta \tau  + 0.0063j{(\delta \tau )^2} + 0.0421j{(\delta \tau )^3}\\
\quad\quad + (0.0362 + 0.0137j){(\delta \tau )^4} + o({(\delta \tau )^4}),
\end{array}
\end{array}} \right. \\
&& \left\{ \begin{array}{l}
\delta \lambda  = (0.0385 + 0.0698j){(\delta \tau )^{\frac{1}{2}}} + o({(\delta \tau )^{\frac{1}{2}}}),\\
\delta \lambda  = 0.1592j\delta \tau  + 0.0253j{(\delta \tau )^2} + 0.6696j{(\delta \tau )^3}\\
\quad\quad  + (1.1585 + 0.4376j){(\delta \tau )^4} + o({(\delta \tau )^4}),
\end{array} \right. \\
&& \left\{ \begin{array}{l}
 \delta \lambda  =  - 0.1592j\delta \tau  + ( - 05371 + 03644j){(\delta \tau )^2} \\
\quad\quad + o({(\Delta \tau )^2}), \\
 \delta \lambda  =  - 0.0988j{(\delta \tau )^{\frac{1}{3}}} + ( - 00356 + 00028j){(\delta \tau )^{\frac{2}{3}}} \\
\quad\quad + o({(\delta \tau )^{\frac{2}{3}}}), \\
 \end{array} \right. \\
&& \left\{ \begin{array}{l}
\delta \lambda  =  - 0.0796j\delta \tau  + ( - 0.0671 +
0.0487j){(\delta \tau )^2} \\
\quad\quad  + o({(\delta \tau )^2}),\\
\delta \lambda  =  - 0.1592j\delta \tau  + 0.0253j{(\delta \tau )^2} + 0.6615j{(\delta \tau )^3}\\
\quad\quad  + ( -1.1585 - 0.4363j){(\delta \tau )^4} + o({(\delta \tau )^4}),
\end{array} \right. 
\end{eqnarray*}
for $k=0, 1, 2,$ and  $3$, respectively.
\begin{figure}
  \centering
  \includegraphics[height=1.7in,width=2.5in]{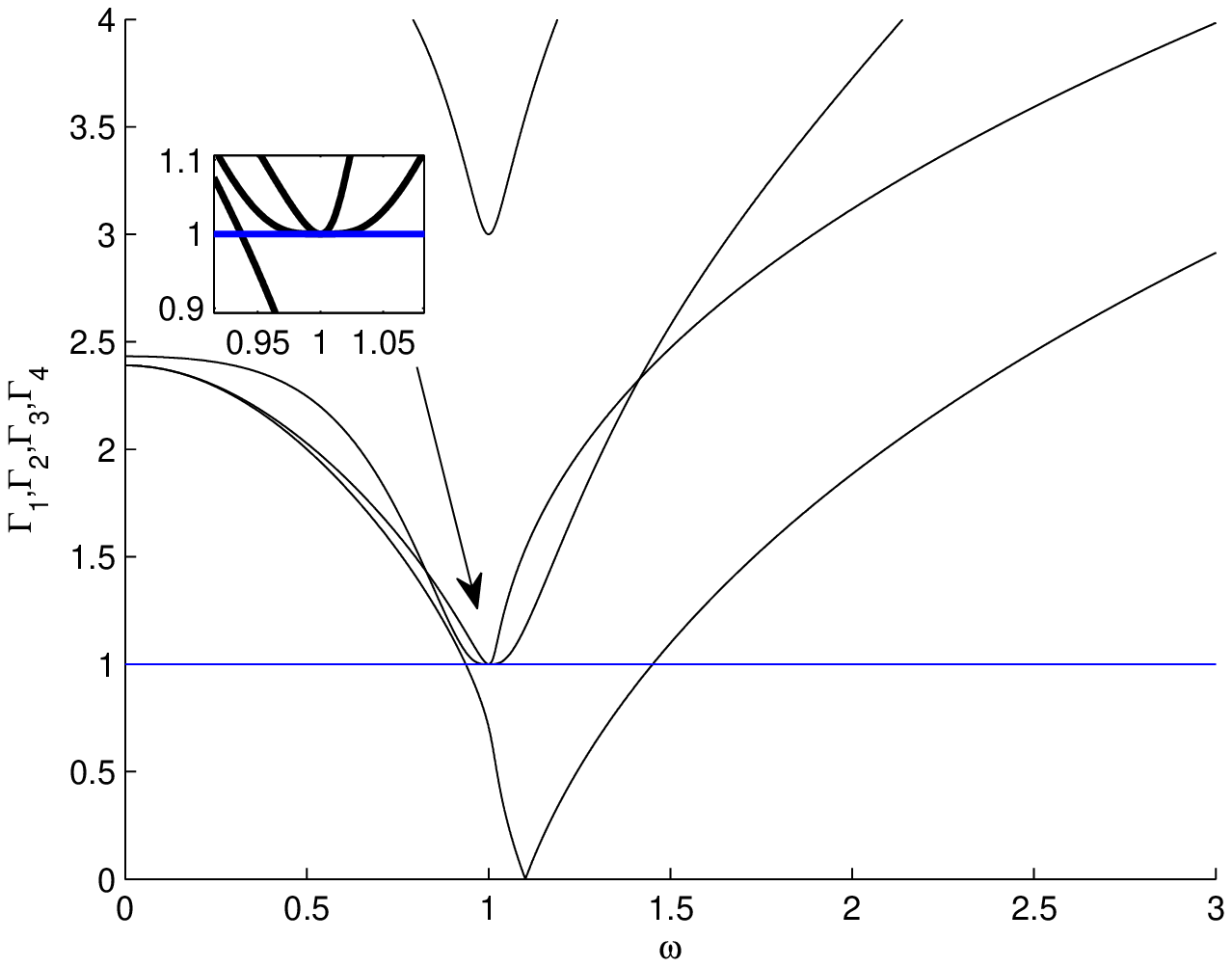}\\
  \caption{Example \ref{MotivatingExampleVariantMultiplicityTACFP2017}: FSCs} \label{FSCExampleInvarianceConfirmationofMotivatingExampleVariantMultiplicityTACFP2017}
\end{figure}
The above Puiseux series are consistent with the analysis by Theorem \ref{InvarianceTheorem}. In this case, for each $k$, (i) the Puiseux series has multiple conjugacy classes; (ii) the Puiseux series involves many degenerate terms, and, finally, (iii) the structure of Puiseux series is variable w.r.t.\ different $k$. It is evident that the FSC formalism significantly simplifies the analysis. $\Box$
\end{example}

In the sequel, we present an example reported in \cite{Li-etal:19} to explicitly illustrate that the asymptotic behavior of a multiple characteristic root may lead to a \emph{stability reversal}.

\begin{example}\label{exdouble}
Consider the quasipolynomial $\Delta(\lambda;\tau)=\lambda^5-\sum_{\ell=0}^{4}\alpha_{\ell}\lambda^{\ell}-(\sum_{\ell=0}^{4}\beta_{\ell}\lambda^{\ell})e^{-\tau\lambda}$.
where $\alpha_{0}=\frac{\pi}{2}-\frac{\pi^{2}}{8}-1$, $\alpha_{1}=-2+\frac{\pi}{2}$, $\alpha_{2}=-\frac{\pi^{2}}{4}+\pi-10$, $\alpha_{3}=-3+\frac{\pi}{2}$, $\alpha_{4}=-\frac{\pi^{2}}{8}+\frac{\pi}{2}-8$, $\beta_{0}=-1$, $\beta_{1}=-1$, $\beta_{2}=-10$, $\beta_{3}=-1$, and $\beta_{4}=-8$.\\
For this  system, ${\rm NU}(\tau)=0$ for $\tau\in[0,1.2525)$ and ${\rm NU}(\tau)=2$ if $\tau\in(1.2525,\pi)$. At $\tau=\pi$, a critical imaginary root $\lambda=j$ with $n=2$ and $g=1$ appears, whose asymptotic behavior corresponds to the  Puiseux series:
\begin{equation*}
\delta\lambda=0.1468j(\delta\tau)^{\frac{1}{2}}+(-0.0033-0.1473j)(\delta\tau)+\ldots
\end{equation*}
It implies that when $\tau$ passes through $\pi$, a stability reversal is caused by the asymptotic behavior of the double root $\lambda=j$. The system possesses two and only two stability intervals of $\tau$. More precisely, it is asymptotically stable iff $\tau\in[0,1.2525)\cup(\pi,4.0549)$.
\begin{figure}[!t]
\centerline{
\includegraphics[height=1.7in,width=1.8in]{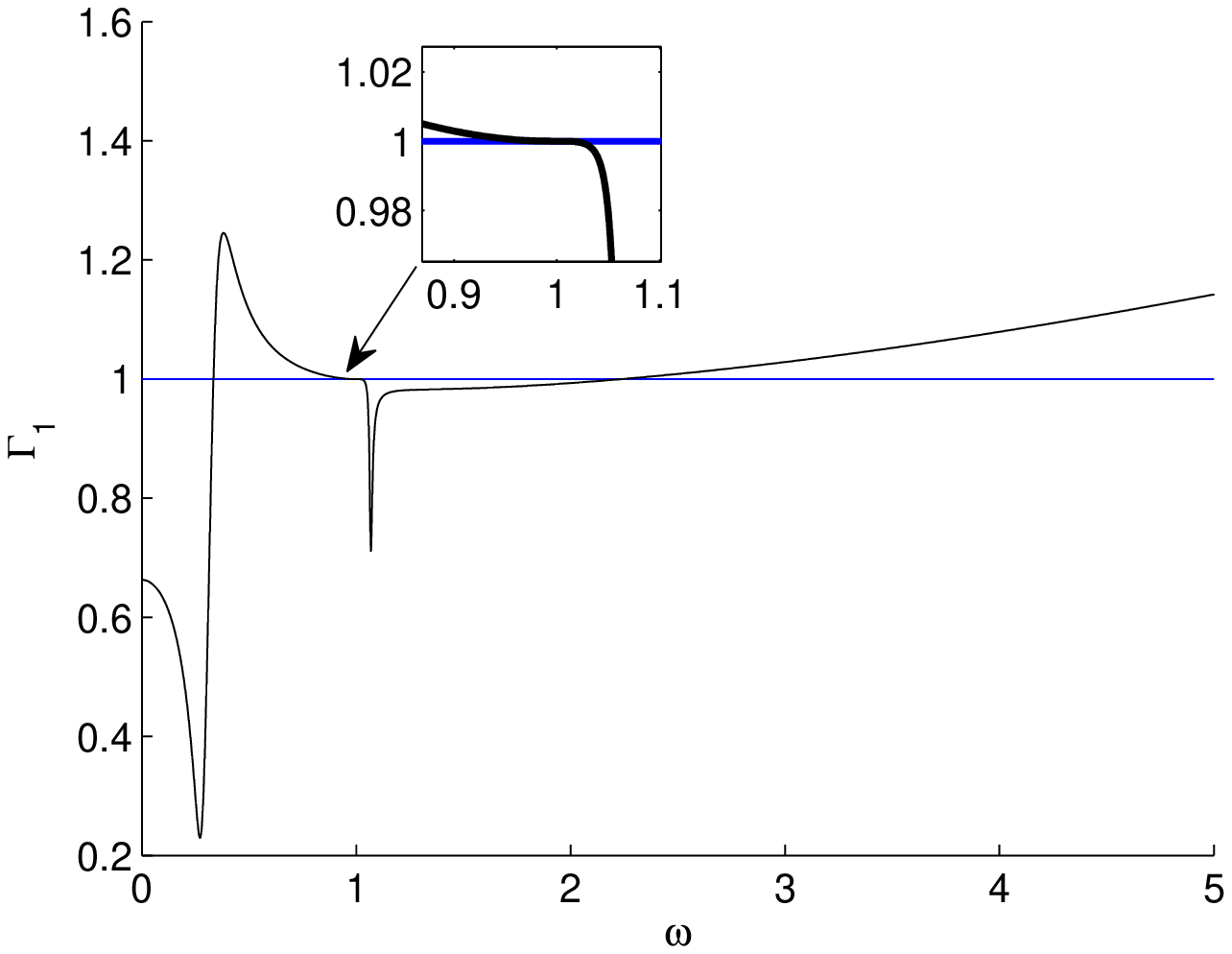}
\hfil \includegraphics[height=1.7in,width=1.8in]{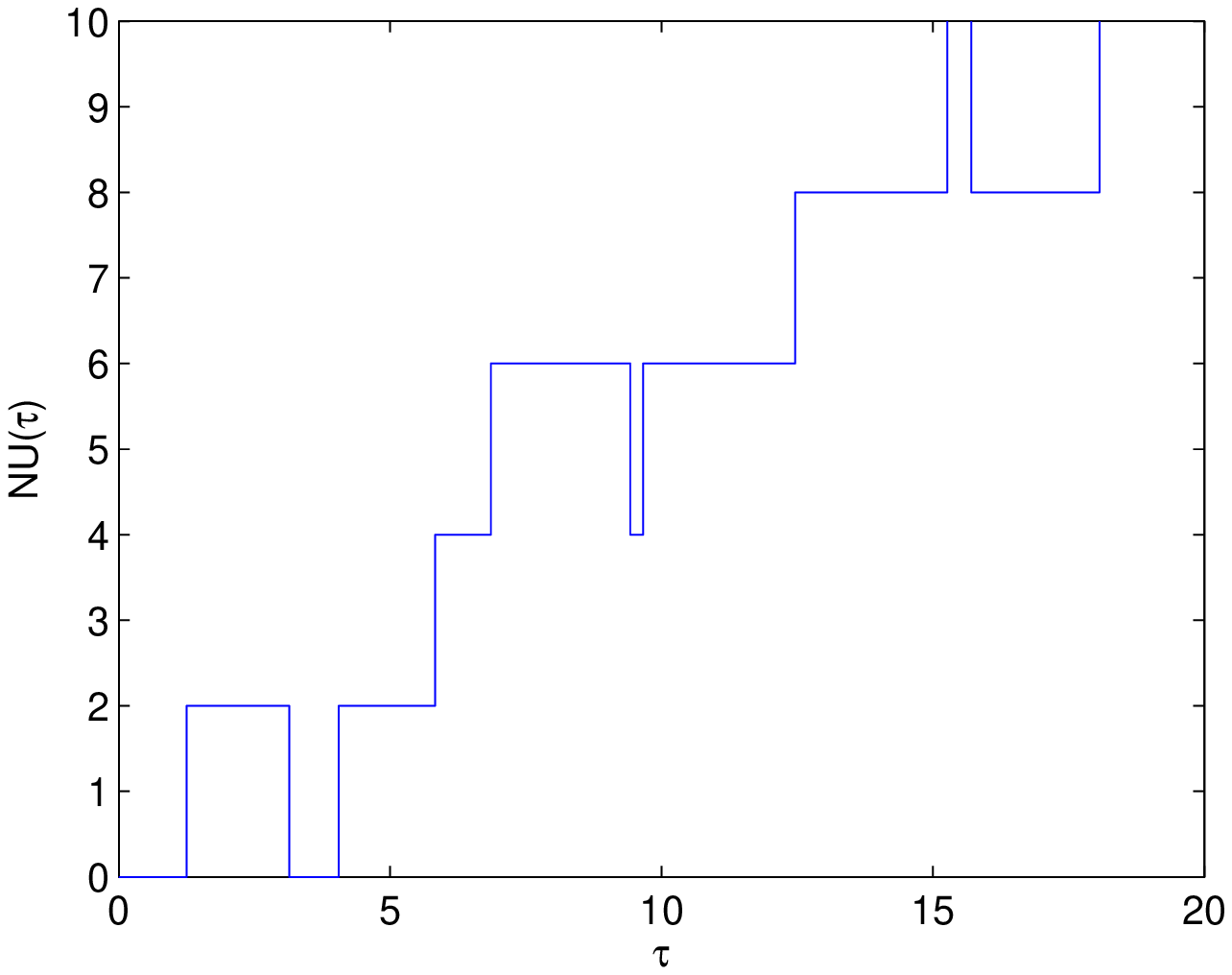}
}
\caption{Example  \ref{exdouble}: FSCs and $NU(\tau)$ }\label{FSCandNUExampleLXAutomatica}
\end{figure}
The asymptotic behavior issue and the stability can be easily addressed by using the frequency-sweeping approach (see Fig.~\ref{FSCandNUExampleLXAutomatica}). $\Box$
\end{example}

\section{Multiple characteristic roots and multiple delays}
\label{sec:mult-del}

Let us consider the following quasipolynomial:
 \begin{align}\label{eq:quasipolynomial}
  \Delta(\lambda,\overrightarrow{\tau})&=\sum_{k=0}^{2}P_k(\lambda)e^{-\tau_k \lambda},\qquad\tau_k\geq0,
 \end{align}
 where $\overrightarrow{\tau}:=\left(\tau_1,\tau_2\right)$,  $\tau_0:=0$, $\tau_1,\tau_2\in\mathbb{R}_+$ and,
  \[
  P_{0}(\lambda):=\lambda^n+\sum_{\ell=0}^{n-1}a_{0\ell}\lambda^{\ell},\,\, P_{k}(\lambda):=\sum_{\ell=0}^{n-1}a_{k\ell}\lambda^{\ell},
  	   \quad k=1,2.
 \]
 This section focuses on computing the first approximation of the solution of quasipolynomials around multiple imaginary roots. To such an end, we will derive conditions allowing to express the solutions of $\Delta$ as a Puiseux series expansion, that is:
              \[
               \lambda(\tau_1,\tau_2)=c(\tau_2^{1/d})\tau_1^{\beta}+o(\tau_1^{1/d}\tau_2^{1/d}), 
              \]
              where $\beta=\alpha/d$ and $\alpha\in\mathbb{N}$.\\
Inspired by the approach adopted in Section~\ref{sec:Pui}, we will compute the associated Weierstrass polynomial of $\Delta$. In this vein, following \cite{Mailybaev:01}  for the case of a holomorphic function, $f(z,\overrightarrow{x})$ of complex variables with 
 $\overrightarrow{x}=(x_1,x_2)$ and $z=0$ a $m$-multiple root at $\left(x_1,x_2\right)=\left(0,0\right)$ the
 computation is given as follows. The coefficients $w_i$ \eqref{eq:localf} are analytic, $w_i(0,0)=0$
 and can expressed as convergent power series: 
 \[
  w_i(x_1,x_2)=\sum_{h_1+h_2=1}^{\infty}\frac{1}{h_1!h_2!}w_{i,\overrightarrow{h}}x_1^{h_1}x_2^{h_2}, 
 \]
 where $\overrightarrow{h}=(h_1,h_2)$. It is not difficult to see that the coefficients $w_{i,\overrightarrow{h}}$
 can be computed as follows:
  \[
  w_{i,\overrightarrow{h}}=\left.\frac{\partial^{h_1+h_2} w_i}{\partial x_2^{h_1}\partial x_1^{h_2}}\right|
    _{(0,0)}.
 \]
 Since the analytic function locally satisfies $f=Wb$, its partial derivatives satisfy the following recursive relations:
 \begin{align}\label{eq:wicoef}
   w_{i,\overrightarrow{h}}&=\sum_{r=0}^{i}\alpha_{ir}F_{r,\overrightarrow{h}},\\
   F_{r,\overrightarrow{h}}&=f_{r,\overrightarrow{h}}-\sum_{k=0}^r\sum_{\overrightarrow{h}^{\prime}+\overrightarrow{h}^{\prime\prime}=\overrightarrow{h}}\!\!\!\!\!c\!\left(r,k;\overrightarrow{h}^{\prime},\overrightarrow{h}^{
                    \prime\prime}\right)w_{k,\overrightarrow{h}^{\prime}}b_{r-k,\overrightarrow{h}^{\prime\prime}}
                    ,\nonumber 
 \end{align}
 with $\overrightarrow{h}^{\prime}\neq \overrightarrow{0}$, $\overrightarrow{h}^{\prime\prime}\neq \overrightarrow{0}$  and constant coefficients:
 \[
  \alpha_{rr}:=\frac{m!}{r!f_{m,\overrightarrow{0}}},\quad \alpha_{ir}:=-\frac{m!}{f_{m,\overrightarrow{0}}}
    \sum_{k=r}^{i-1}\frac{f_{m+i-k,\overrightarrow{0}}\alpha_{kr}}{(m+i-k)!},
 \]
 \[
  c\left(r,k;\overrightarrow{h}_1,\overrightarrow{h}_2\right):=\frac{r!}{(r-k)!}\prod_{s=1}^2 
					      \frac{(h_s^{\prime}+h_s^{\prime\prime})!}{h_s^{\prime}!h_s^{\prime\prime}!},\\
 \]
 and for $\overrightarrow{h}^{\prime}\neq0$, $k^{\prime}=k+m$, $b_{k,\overrightarrow{h}}$ is given by
 \[
  \!\!\frac{k!}{(m\!+\!k)!}\!\!\!\left[\!f_{k^{\prime},\overrightarrow{h}}-\!\!\sum_{r=0}^{m-1}\sum_{\overrightarrow{h}^{\prime}+\overrightarrow{h}
                    ^{\prime\prime}=\overrightarrow{h}}\!\!\!\!\!\!\!c\!\left(k^{\prime},r;\overrightarrow{h}^{\prime},\overrightarrow{h}^{\prime\prime}\!\right)\!\!w_{r,
                    \overrightarrow{h}^{\prime}}b_{k^{\prime}-j,\overrightarrow{h}^{\prime\prime}}\!\right]\!\!.
 \]
 In the sequel, we will adopt the following notations.  
\begin{definition}\label{def:ni}
  Let the natural numbers $n_i^{(r)}$, for $i\in\llbracket0,m-1\rrbracket$ and $r=1,2$, denote the first nonzero  partial derivative in $(z,x_1,x_2)$ of $f$, such that:
\[
   f(0,0,0)\!=\!\frac{\partial^{i}f}{\partial z^i}\!=\cdots=\!\frac{\partial^{i+n_i^{(r)}-1}f}{\partial 
        z^i\partial\tau_r^{n_i^{(r)}-1}}=0, \quad 
        \frac{\partial^{i+n_i^{(r)}}f}{\partial z^i\partial\tau_j^{n_i^{(r)}}}\neq0,
  \]
  with derivatives evaluated at $(0,0,0)$. For $n_i^{(r)}=\infty$ we have derivatives
  \[
   \frac{\partial^{i}f}{\partial z^i}\!=\cdots=\frac{\partial^{i+n_i^{\prime}-1}f}{\partial z^i
    \partial\tau_2^{n_i^{\prime}-1}}=0, \quad \frac{\partial^{i+n_i^{\prime}}f}{\partial z^i
    \partial\tau_2^{n_i^{\prime}}}\neq0,    
  \]
  evaluated at $(z,\overrightarrow{x})=(0,0,1)$, with $n_i^{\prime}\in\mathbb{N}\cup\{0\}$.
\end{definition}

  Leading terms of coefficients $w_i$ can be easy found up to the $n_i^{(r)}$ and $n_i^{\prime}$ derivatives. More precisely, as a first observation we give the following result: 
 \begin{proposition}\label{prop:lterms}
  Suppose that the Weierstrass polynomial has first nonzero partial derivative, such that 
  \[
   n_i^{(r)}>n_{i+1}^{(r)}, \quad 0\leq i<m \text{ and } r=1,2.
  \]
  Then, the leading terms of $w_i(\overrightarrow{x})$ are given by
  \[
   w_i(x_1,x_2)=\alpha_{i,i}f_{i,(n_i^{(1)},0)}x_1^{n_i^{(1)}}+\alpha_{i,i}f_{i,(0,n_i^{(2)})}
               x_2^{n_i^{(2)}}+\cdots.
  \] 
  If $n_i^{(r)}=\infty$, we have
  \[
   w_i(x_1,x_2)=\alpha_{i,i}f_{i,(n_i^{\prime},\eta)}x_1^{n_i^{\prime}}x_2^{\eta}+\cdots.
  \]   
 \end{proposition}
  \begin{remark}\label{rem:invroots}
 As in the single parameter case, there may be a situation in which 
  \[
   \left.f_{i,\overrightarrow{h}}\right|_{(0,0,0)}=0 \quad \forall \;h_1,h_2\in\mathbb{N}.
  \]
  Since $w_i$ are analytic functions, this is equivalent to $w_i(\overrightarrow{x})\equiv0$ for $i\leq i\leq
  \kappa-1$. Thus, according to Theorem \ref{theo:Weierstrass} $f$ has the following local structure:
  \[
   z^{\kappa}\left[z^{m-\kappa}+w_{m-\kappa}(\overrightarrow{x})z^{m-\kappa-1}+\cdots+w_{\kappa}(\overrightarrow{x})
                   \right]b(z,\overrightarrow{x}).                   
  \]
  Thus, there are $\kappa$-invariant solutions $z=0$ for all $\overrightarrow{x}$. If such number $\kappa$ does not exist
  (i.e., if such situation does not occur), then $\kappa$ will be simply defined as $\kappa:=0$. $\Box$
 \end{remark}
 
 \subsection{Newton diagram method: An extension}
 
 For a given $\pi_k\in\Pi$, we consider the order of $w_k$ in
 $x_1$, taking $x_2$ as an element of $\C\left[\left[x_2^{1/d}\right]\right]$ (for an appropriate $d\in\mathbb{N}$). For such a purpose, the following
 definition will be useful:
  \begin{equation}\label{eq:order}
  \rho_k:=\ord{w_k(x_1,x_2)}{x_1}=\ord{w_k(x_1,1)}{}.
 \end{equation}
 Then, the \emph{Newton polygon} of $W_p(z,\overrightarrow{x})$, with respect to $x_1$, is defined by the lower boundary of the convex hull of the points $(k,\rho_k)\in\Pi$. Hence, in order to apply the the Newton diagram procedure, the solution $z$ will take the following structure
 \[
  z(x_1,x_2)=\sum_{i}c_i(x_2)x_1^{i/d},
 \]
 where the coefficient $c_i(x_2)$, is in general, given by a single parameter Puiseux series in $x_2$.\\
 Therefore, with the aim to compute such a Puiseux series, let us suppose that we have determined the Newton diagram of the Weierstrass polynomial of $\Delta$. Since we are dealing with a monic polynomial, the Newton polygon has a finite number of segments,
 each one with a corresponding set of points $\Pi^{(\ell)}$ and rational numbers $\beta_{\ell}\geq0$ satisfying
 \[
  \beta_0>\beta_1>\cdots >\beta_{r}.
 \]
 In this regard, the segments are presented in two possible ways. The first one corresponds to a Newton polygon with a
 horizontal segment with $\beta_i=0$, and the second one where $\beta_j>0$ (for $i\neq j$). In this vein, for
 $0\leq\ell<m$, the Newton Diagram $\Pi$ is given as the set $\Pi=\Pi^{\prime}\cup\!\Pi^{\prime\prime}$:
 \[
  \left\{\!(0,\rho_0),\dots,\left(\ell,0\right)\!\right\}\cup
                                        \left\{\!\left(\ell,0\right),\dots,(k,\rho_{k}),\dots,
                                        \left(m,0\right)\!\right\}.
 \]
 Assuming that at the first step of the process we found a horizontal segment with a slope $\beta_r=0$, the following result allows characterizing the solution structure.
  \begin{proposition}\label{prop:hsegment}
  Let $W_p(z,\overrightarrow{x})$ be a Weierstrass polynomial of a given quasipolynomial $\Delta$. Suppose that at least one coefficient $w_i(\overrightarrow{\tau})$ possesses order $\rho_i=0$. Then, the equation $\mathcal{P}_i(z;\tau_1)=0$ of the corresponding $i$-horizontal segment has solutions
  $c_k(\tau_2^{1/d})$ in the form of Puiseux series.
 \end{proposition}
 \begin{remark}
  In the above result, $\mathcal{P}_i$ is built in the same fashion as in \eqref{eq:pnu}. $\Box$
 \end{remark}
 Now, assuming that at the first step of the process we found a segment with a negative slope $\beta_r<0$, the following result proposes an appropriate change of variable that will allow us to pursue the procedure.
  \begin{proposition}\label{prop:nslope}
  Let $W_p(z,\overrightarrow{\tau})$ be a Weierstrass polynomial of a given quasipolynomial $\Delta$ and assume that the first Newton diagram possesses a segment with negative slope. Then, there exist a change of variables 
  $(\lambda,\tau_1,\tau_2)\mapsto (\xi,y_1,y_2)$ such that the polynomial $\mathcal{P}_i(z;y_2)$ has Puiseux series solutions $c_k(y_2^{1/d})$. 
 \end{proposition}
 
 \subsection{Puiseux series expansion}
 
 Since any critical solution $\left(\lambda^{\ast},\tau_1^{\ast},\tau_2^{\ast}\right)$ can always be translated to the origin by appropriate shifts $\lambda\mapsto \lambda-\lambda^{\ast}$, $\tau_1\mapsto\tau_1-\tau_1^{\ast}$, $\tau_2
 \mapsto\tau_2-\tau_2^{\ast}$, hereinafter we will assume that $\left(\lambda^{\ast},\tau_1^{\ast},\tau_2^{\ast}\right)=(0,0,0)$.
 \begin{proposition}\label{prop:asymp}
  Let $\lambda^{\ast}=j\omega^{\ast}$ be a $m$-multiple root of 
  $\Delta(\lambda,\overrightarrow{\tau})$ at $\overrightarrow{\tau}^{\ast}:=\left(\tau_1^{\ast},\tau_2^{\ast}\right)$. Assume that $\kappa=0$ and $r$, $\beta_h$, $\left(i_h,\ell_h\right)$, $m_h$ and $\Pi^{(h)}$, for $h\in\llbracket0,r-1\rrbracket$ are given by the Algorithm \ref{alg}. 
  Then, at $\overrightarrow{\tau}=\overrightarrow{\tau}^{\ast}$ the $m$-zeros of $\Delta(\lambda,\overrightarrow{\tau})$ can be expanded as 
 \begin{multline*}
  \lambda_{hq}(\overrightarrow{\tau}) = 
  j\omega^{\ast}+c_{hq}(\tau_2) \left(\tau_1 - \tau_1^{\ast}\right)^{\beta_{h}} \\
	 +o \left( \left|\tau_1 - \tau_1^{\ast}\right|^{\beta_{h}} \left|\tau_2 - \tau_2^{\ast}\right|^{\beta_{h}^{\prime}}\right),
 \end{multline*}
 for $h\in\llbracket0,r-1\rrbracket$, $q\in\llbracket0,m_{h}\rrbracket$ and $m=m_0+\cdots+m_{r-1}$. For $\beta_h>0$, $c_{hq}(\tau_2)$ are 
 roots of the polynomial:
 \begin{equation*}
  \mathcal{P}_{h}(z;\tau_2)\!=\!\!\!\sum_{k=i_{h-1}}^{i_{h}}\!\!\!w_{k,(n_k^{(1)},n_k^{\prime})}\tau_2^{n_k^{\prime}}z^{k-i_{h-1}}\!,
	  \quad\left(k,n_k^{(1)}\right)\!\in\!\Pi^{(h)},
 \end{equation*}
 when $\beta_{r-1}=0$, the coefficients are given by the solution of
 \begin{equation*}
  \mathcal{P}_{h}(z;\tau_2)\!=\!\!\!\sum_{k=i_{h-1}}^{i_{h}}\!\!\!w_{k,(0,n_k^{(2)})}\tau_2^{n_k^{(2)}}\!z^{k-i_{q-1}}\!,\quad\left(k,0\right)\!\in\!\Pi^{(r-1)},
 \end{equation*}
 where $n_k^{(1)}$, $n_k^{(2)}$, $n_k^{\prime}$ are given by the first nonzero partial derivatives of Definition \ref{def:ni}; the constant terms $w_{k,(n,\eta)}\in\C$ are computed using \eqref{eq:wicoef}.
 \end{proposition}
 To illustrate the previous results, consider the following:
 \begin{example}
 Consider the quasipolynomial $\Delta(\lambda;\tau_1,\tau_2)\allowbreak:=P_0(\lambda)+P_1(\lambda)e^{-\tau_1\lambda}+P_2(\lambda)e^{-\tau_2\lambda}$
 where
 \begin{subequations}\label{eq:exp2}
 \begin{eqnarray}
  P_0(\lambda)&:=&\lambda^5\!+\!\lambda^4\!+\!\frac{4+\pi}{2}\lambda^3\!+\!2\lambda^2\!+\!\frac{2+\pi}{2}\lambda\!+\!2,\\
  P_1(\lambda)&:=&1, \quad P_2(\lambda):=2\lambda^4\!+\!4\lambda^2+2.
 \end{eqnarray}
 \end{subequations}
 For $\overrightarrow{\tau}^{\ast}=(\pi,1)$, $\Delta$ has a double root at $\lambda=j$. The shifting from $(j,\pi,1)$ to the origin leads to $\widetilde{\Delta}$. Next, by computing the first nonzero partial derivatives of the quasipolynomial at $(0,0,0)$, one gets:
 \[
  \Delta_{0,(1,0)}=j\;\,\,\,\,\,\,\,\,\,\,\,\,\,\,\Rightarrow n_0^{(1)}=1,\quad \Delta_{0,(0,n)}=0 \Rightarrow
                       n_0^{(2)}=\infty,
 \]
 \[
  \Delta_{1,(1,0)}=1-j\pi \Rightarrow n_1^{(1)}=1, \quad \Delta_{1,(0,n)}=0 \Rightarrow n_1^{(2)}=\infty.
 \]
 Hence, by Proposition \ref{prop:lterms}, we have that  $\rho_h=n_h^{(1)}$ for $h=1,2$. The Newton
 diagram is given by $\Pi=\{(0,1),(1,1),(2,0)\}$. Table \ref{tab:exp2a} summarizes the results deriving from Algorithm \ref{alg}.
 \begin{table}[h!]
 \centering
 \caption{Results summary for $\Delta$ given by \eqref{eq:exp2}.}
  \scalebox{0.79}{\begin{tabular}{l l l} 
  \toprule
   Initial Data & Algorithm Output & $\mathcal{Z}:=\left\{z\in\C:\;\mathcal{P}_h(z)=0\right\}$\\ 
  \midrule
   $m=2$, $ \kappa=0$. $\rho_0=1$ & $r=1$, $m_{0}=2$, $\beta_{0}=1/2$ & $\mathcal{P}_0(z):=z^2+w_{0,(1,0)}$\\
   $\Pi=\left\{(0,1),(1,1),(2,0)\right\}$ & $\Pi^{^{(0)}}=\left\{(0,1),(2,0)\right\}$ & $\left\{c_{0,q}=\pm\sqrt{w_{0,(1,0)}}\right\}$\\
  \bottomrule
  \end{tabular}}\label{tab:exp2a}
 \end{table}
 
 Since $w_1$ is not over the Newton polygon, following Proposition \ref{prop:lterms}, we compute the leading term  of $w_i(\overrightarrow{\tau})$ as:
 \[
  w_0(\overrightarrow{\tau})=\frac{-2j}{(8+\pi^2)+j(8-3\pi)+16e^{-j}}\tau_1+\cdots.
 \]
 From the algorithm output, we get a segment with slope $\beta_0=1/2$. According to Proposition \ref{prop:asymp}, we need to solve:
 \[
  \mathcal{P}_0(z)=z^2-\frac{2j}{(8+\pi^2)+j(8-3\pi)+16e^{-j}}=0,
 \]
and, for $q\in\{0,1\}$, the solutions are given by
 \begin{equation*}
  \lambda_{0,q}(\overrightarrow{\tau})\!=\!j+\frac{(-1)^{q}
  \sqrt{2}j^{3/2}}{\sqrt{(8\!+\!\pi^2)\!+\!j(8\!-\!3\pi)\!+\!16e^{-j}}}(\tau_1\!-\!\pi)^{1/2}+\cdots.
 \end{equation*}
\end{example}

\begin{remark}
It is worth mentioning that the frequency-sweeping approach can be also extended to handle multiple characteristic roots on imaginary axis for delay systems including incommensurate delays, see, e.g. \cite{lnc:19-ima} for an \emph{iterative frequency-sweeping method}. $\Box$
\end{remark}

\section{Hypergeometric functions and multiplicity-induced-dominancy}
\label{sec:MID}

We consider in this section the DDE given by:
\begin{equation}
\label{eq:general-DDE-single-delay}
y^{(n)}(t) + \sum_{k=0}^{n-1} a_k y^{(k)}(t) + \sum_{k=0}^m \alpha_k y^{(k)}(t - \tau) = 0,
\end{equation}
where $y$ is real-valued, $n$ is a positive integer, $m \in \llbracket 0, n\rrbracket$, $\tau > 0$ is the delay, and $a_0, \dotsc, a_{n-1}, \alpha_0, \dotsc, \alpha_m$ are real coefficients. Its characteristic function is
\begin{equation}
\label{eq:Delta-DDE-single-delay}
\Delta(\lambda;\overrightarrow{a},\overrightarrow{\alpha}) = \lambda^n + \sum_{k=0}^{n-1} a_k \lambda^k + e^{-\lambda \tau} \sum_{k=0}^m \alpha_k \lambda^k.
\end{equation}
The problem of characterizing regions in the space of parameters of \eqref{eq:general-DDE-single-delay} ensuring exponential stability is a highly non-trivial problem of ongoing interest, and, thanks to classical results, such a problem is equivalent to characterizing regions in the space of the coefficients of \eqref{eq:Delta-DDE-single-delay} ensuring that all its roots $\lambda$ satisfy $\Re(\lambda) \leq -\gamma$ for some $\gamma > 0$ (see, e.g., \cite{mn14:siam}). This question is naturally related to control-theoretical problems, since \eqref{eq:general-DDE-single-delay} can be seen as the closed-loop system obtained by applying a linear (possibly delayed) feedback law to a controlled delay-differential equation, and in this case, by suitably choosing the coefficients of the feedback law, one may choose (some of) the coefficients of \eqref{eq:general-DDE-single-delay}.\\
Since \eqref{eq:Delta-DDE-single-delay} admits infinitely many roots but has only $m + n + 1$ parameters, one cannot expect to be able to choose arbitrarily the location in $\mathbb C$ of all roots of $\Delta$. Some works are interested in methods to choose the location of finitely many roots of $\Delta$, trying to guarantee that the other roots have negative real part and are separated from the imaginary axis, a technique known as \emph{partial pole placement}. Trial-and-error methods, such as those in \cite{Ram2011Partial}, have been proposed in the literature, but techniques guaranteeing that the non-assigned roots have negative real part are more difficult to obtain. The numerical paradigm known as \emph{continuous pole placement}, introduced in \cite{Michiels2002Continuous}, exploits continuity of roots of $\Delta$ with respect to the coefficients of the system in order to move roots with positive real part to the left half-plane while ensuring that no stable roots becomes unstable. 

In the sequel, we present a recent technique for partial pole placement, based on the property known as \emph{multiplicity-induced-dominancy}, or MID for short (see, e.g., \cite{Boussaada2020Multiplicity}). The MID property states that, if a real root of \eqref{eq:Delta-DDE-single-delay} attains its maximal multiplicity (which is equal to $m + n + 1$, cf.\ Remark~\ref{remk:PolyaSzego}), then this root necessarily is the rightmost root in the complex plane. Hence, a technique for partial pole placement based on the MID property consists in choosing the coefficients of the system in order to ensure the existence of a negative real root of maximal multiplicity. More precisely, we have the following:
\begin{theorem}
\label{thm:MID}
Consider the quasipolynomial $\Delta$ given by \eqref{eq:Delta-DDE-single-delay} and let $\lambda_0 \in \mathbb R$. The number $\lambda_0$ is a root of maximal multiplicity $m+n+1$ of $\Delta$ if and only if
\begin{equation}
\label{Coeffs}
\begin{aligned}
a_k & = \left( -1 \right) ^{n-k}n!\,\sum_{i=k}^{n}{\frac {\binom{i}{k} \binom{m+n
-i}{m} \lambda_0^{i-k}}{i!\,{\tau}^{n-i}}}, 
\quad k \leq n-1,\\
\alpha_k & = \left( 
-1 \right) ^{n-1}{{e}^{\lambda_{{0}}\tau}}\sum_{i=k}^{m}{\frac {
 \left( -1 \right) ^{i-k} \left( m+n-i \right) !\,\lambda_0^{i-k
}}{k!\, \left( i-k \right) !\, \left( m-i \right) !\,{\tau}^{n-i}}}, 
k \leq m.
\end{aligned}
\end{equation}
If $m < n$ and \eqref{Coeffs} is satisfied, then $\Re(\lambda) < \lambda_0$ for every root $\lambda \neq \lambda_0$ of $\Delta$. If $m = n$ and \eqref{Coeffs} is satisfied, then $\Re(\lambda) = \lambda_0$ for every root $\lambda$ of $\Delta$. In particular, the trivial solution of \eqref{eq:general-DDE-single-delay} is exponentially stable if and only if $a_{n-1} > -\frac{n (m+1)}{\tau}$.
\end{theorem}

Theorem~\ref{thm:MID} was proved in the case $m = n - 1$ in \cite{MBN-2021-JDE} and extended to any $m \in \llbracket 0, n\rrbracket$ (including thus the neutral case $m = n$) in the recent paper \cite{BoussaadaGeneric}. Let us briefly present the main ideas of its proof. Up to a change of variables corresponding to a translation and a scaling, it suffices to prove the theorem in the case $\tau = 1$ and the desired root of maximal multiplicity is at the origin, i.e., $\lambda_0 = 0$. In this case, $\lambda_0 = 0$ is a root of maximal multiplicity $m + n + 1$ if and only if $\Delta(0) = \dotsb = \Delta^{(m + n)}(0) = 0$, which gives a linear system in the $m + n + 1$ coefficients $a_0, \dotsc, a_{n-1}, \alpha_0, \dotsc, \alpha_m$ admitting \eqref{Coeffs} as its unique solution.

To prove the properties on the dominance of $\lambda_0$, the main ingredient is that, under \eqref{Coeffs} and with $\tau = 1$ and $\lambda_0 = 0$, $\Delta$ can be factorized as
\[
\Delta(\lambda;\overrightarrow{a},\overrightarrow{\alpha}) = \frac{n! \lambda^{m + n + 1}}{(m + n + 1)!} \Phi(m + 1, m + n + 2, -\lambda),
\]
where $\Phi$ is Kummer confluent hypergeometric function, which admits the integral representation
\[
\Phi(a, b, z) = \frac{\Gamma(b)}{\Gamma(a) \Gamma(b - a)} \int_0^1 t^{a - 1} (1 - t)^{b - a - 1} e^{z t} d t
\]
for every $a, b, z \in \mathbb C$ with $\Re(b) > \Re(a) > 0$, where $\Gamma$ denotes the Gamma function. In particular, proving the last part of Theorem~\ref{thm:MID} amounts to showing that all roots of $\Phi(m + 1, m + n + 2, \cdot)$ have positive real part if $m < n$ and real part zero if $m = n$. As detailed in \cite{BoussaadaSome}, these facts can be established by applying the technique developed in \cite{Hille1922Oscillation} to a family of special functions related to Kummer functions, known as Whittaker functions, and exploring the fact that Whittaker functions satisfy a second-order ordinary differential equation.

\subsection{Inverted pendulum: Delay, PD control \& MID}

The MID property may hold even when the multiplicity of a given quasipolynomial's root does not reach its maximal value. Hereafter, we illustrate such a claim on the standard comprehensive problem of the stabilization of the inverted pendulum. The equation of motion of an inverted pendulum controlled by a delayed PD-controller writes as:
\begin{equation} \label{eq_mot_sing_inv_pend}
\begin{aligned}
&\ddot \varphi(t)+a_0\varphi(t)=u(t)\,,\\
&u(t)=-b_0 \varphi(t-\tau)-b_1 \dot\varphi(t-\tau)\,,
\end{aligned}
\end{equation}
with $\tau>0$ and $a_0<0$. The characteristic function corresponding to \eqref{eq_mot_sing_inv_pend} is:
\begin{equation} \label{char_eq_sing_inv_pend}
\Delta(\lambda;b_0,b_1,\tau)=P_0(\lambda)+P_1(\lambda)\,e^{-\lambda\tau},
\end{equation}
with $P_0(\lambda)=\lambda^2+a_0$ and $P_1(\lambda;b_0,b_1)=b_0+b_1 \lambda$. 
The critical delay of system \eqref{eq_mot_sing_inv_pend} is given by:
\begin{equation} \label{crit_delay_sing_inv_pend}
\tau_\mathrm{crit}=\sqrt{-\tfrac{2}{a_0}}\,,
\end{equation}
that is, the closed-loop system \eqref{eq_mot_sing_inv_pend} is asymptotically stable if and only if $\tau<\tau_\mathrm{crit}$. The critical delay \eqref{crit_delay_sing_inv_pend} can be obtained by studying the multiple roots of $\Delta(\cdot;b_0,b_1,\tau)$. Indeed, assume that $\Delta$ has a real root $\lambda_0$ with algebraic multiplicity at least $\mathrm{deg}P(\lambda)+1=3$. Then $\Delta(\lambda_0;b_0,b_1,\tau)=\Delta'(\lambda_0;b_0,b_1,\tau)= \Delta''(\lambda_0;b_0,b_1,\tau)=0$ give
\begin{equation} \label{sys_sing_inv_pend}
\left\{\begin{aligned}
\lambda_0^2+a_0+e^{-\lambda_0\tau}(b_0+b_1 \lambda_0)&=0\,,\\
2 \lambda_0+e^{-\lambda_0\tau}(-\tau(b_0+b_1 \lambda_0)+b_1)&=0\,,\\
2+e^{-\lambda_0\tau}(\tau^2 (b_0+b_1 \lambda_0)-2\tau b_1)&=0\,.
\end{aligned}\right.
\end{equation}
From \eqref{sys_sing_inv_pend} we obtain
\begin{equation} \label{sys_sing_inv_pend_sol}
\left\{\begin{aligned}
b_0&=e^{\lambda_0 \tau}\left(\tau \lambda_0^3+\lambda_0^2+a_0 \tau \lambda_0-a_0\right) \,, \\
b_1&=-e^{\lambda_0 \tau}\left(\tau \lambda_0^2+2 \lambda_0+a_0 \tau\right)\,, \\
\lambda_0&=\frac{-2\pm\sqrt{2-a_0 \tau^2}}{\tau}\eqqcolon \lambda_{\pm}.
\end{aligned}\right.
\end{equation}
It can be shown that the triple root $\lambda_{+}$ is negative and dominant for every $0<\tau<\tau_\mathrm{crit}$, and therefore system \eqref{eq_mot_sing_inv_pend} is asymptotically stable, see, e.g., \cite{BTNUV-LAA-2017,Boussaada2020Multiplicity}. In particular, at the upper bound $\tau=\tau_{\mathrm{crit}}$ the triple root is $\lambda_{+}=0$ and it is the dominant (rightmost) root of \eqref{char_eq_sing_inv_pend} with control ``gains'' $b_0=-a_0$ and $b_1=-a_0\tau_{\mathrm{crit}}$. 
Alternatively, for a given $\lambda_+=\gamma<0$, \eqref{sys_sing_inv_pend} can be solved for $b_0$, $b_1$ and $\tau$. 
The smallest positive solution for $\tau$ is the critical delay $\tau_{\mathrm {crit}}(\gamma)$ associated with $\gamma$-stability.

\begin{remark}
The dominancy of $\lambda_{+}$ has been shown by using the argument principle, see, for instance \citep{Boussaada2020Multiplicity}.
A constructive proof of the dominancy of $\lambda_{+}$ may be also shown using the corresponding quasipolynomial factorization introduced in \cite{Bedouhene2020Real} and recently extended in \cite{BBIN:22} to arbitrary order systems such that in open-loop they admit only  real-rooted modes. $\Box$
\end{remark}

\subsection{A control oriented MID property}

In this section we consider delayed feedback systems whose characteristic function is a quasipolynomial of the form
\begin{equation} \label{char_eq}
\Delta(\lambda; \overrightarrow{b}, \tau)=P_0(\lambda)+P_1(\lambda,\overrightarrow{b})\,e^{-\lambda\tau},
\end{equation}
where $\tau\in\mathbb{R}_+$, $\mathrm{deg}(P_0)=n$, $\mathrm{deg}(P_1)=n-1$. Assume that the coefficients of $P_0$ are known and that $\overrightarrow{b}\in\mathbb{R}^{n}$ denotes the parameter vector including the coefficients of $P_1$. Assume further that the coefficients $b_i$ are independently adjustable) control parameters:
\begin{equation} \label{pol}
\begin{aligned}
P_0(\lambda)&=a_n \lambda^n+a_{n-1}\lambda^{n-1}+\dotsb+a_1 \lambda+a_0\,, \\[3pt]
P_1(\lambda;\overrightarrow{b})&=b_{n-1}\lambda^{n-1}+b_{n-2}\lambda^{n-2}+\dotsb+b_1 \lambda+b_0\,. \\[3pt]
\end{aligned}
\end{equation}
The problem to be considered is
\emph{to find the values of $\tau$ such that \eqref{char_eq} is $\gamma$-stabilizable.} 
To give a sufficient condition for $\gamma$-stabilizability\footnote{In other words, the spectral abscissa $\alpha_s$ of the closed-loop system should verify the condition $\alpha_s\leq \gamma$ for the corresponding set of parameters.} we utilize the MID-property: the control parameters $b_i$ are tuned such that the characteristic function $\Delta(\lambda)$ has a real root $\lambda_0$ with multiplicity $n+1$.
The result from \cite{BBIN:22} emphasizes the way to factorize a quasipolynomial admitting a multiple root.
\begin{proposition}  \label{Prop__fac}
If the quasipolynomial \eqref{char_eq}-\eqref{pol}
admits a real root $\lambda_0$ with multiplicity at least $n$ then $\Delta(\lambda;\overrightarrow{b},\tau)=$
\begin{equation} \label{fac}
(\lambda-\lambda_0)^n\left(a_n +\int_0^1 e^{-(\lambda-\lambda_0)\tau t}\frac{\tau  R_{n-1}(\lambda_0;\tau t)}{(n-1)!}\mathrm d t\right)\,,
\end{equation}
where the family of polynomials $R_{k}(\lambda;\tau)$ is defined as
\begin{equation} \label{Rpol}
R_k(\lambda;\tau)=\sum_{i=0}^{k} \binom{k}{i} P_0^{(i)}(\lambda)\tau^{k-i},\, k\in\mathbb{N}^{*}\,.
\end{equation}
\end{proposition}
Some sufficient conditions for the dominancy are given by:
\begin{proposition}  \label{Prop__suff_cond}
Let $\lambda_0$ be a real root of the quasipolynomial \eqref{char_eq} with multiplicity at least $n+1$. If $R_{n-1}(\lambda_0;\tau t)\leq 0$, $\forall t$, $0<t\leq1$  then $\lambda_0$ is the dominant root of \eqref{char_eq}.
\end{proposition}

\subsection{Imaginary roots, MID \& applications}

We consider in this subsection a retarded delay-differential equation of second order under the form $\ddot y(t) + a_1 \dot y(t) + a_0 y(t) + \alpha_1 \dot y(t - \tau) + \alpha_0 y(t - \tau) = 0$, where $a_1, a_0, \alpha_1, \alpha_0$ are real parameters to be tuned. Its characteristic function is the quasipolynomial $\Delta$ defined by
\begin{equation}
\label{eq:Delta-mid-complex}
\Delta(\lambda) = \lambda^2 + a_1 \lambda + a_0 +(\alpha_1 \lambda + \alpha_0)\, e^{- \lambda \tau}.
\end{equation}
All previous results on multiplicity-induced-dominancy concern only \emph{real} roots with high multiplicity, and a natural question is whether one may choose as dominant roots a pair of complex-conjugate nonreal multiple roots. The answer, at least for \eqref{eq:Delta-mid-complex}, is affirmative, as stated in the next result.

\begin{theorem}
\label{thm:mid-complex}
Let $\lambda_0 \in \mathbb C$, denote $\sigma_0 = \Re(\lambda_0)$ and $\theta_0 = \Im(\lambda_0)$, and assume that $\theta_0 \neq 0$. Then $\lambda_0$ and its conjugate $\overline \lambda_0$ are roots of multiplicity $2$ of $\Delta$ from \eqref{eq:Delta-mid-complex} if and only if
\begin{subequations}
\label{CoeffsComplex}
\begin{align}
a_1 & \textstyle = - 2 \sigma_{0} - \theta_{0} \frac{2 \tau \theta_{0} - \sin(2 \tau \theta_{0})}{\tau^{2} \theta_{0}^{2} - \sin^{2}(\tau \theta_{0} )}, \displaybreak[0] \label{CoeffsComplexA1}\\
a_0 & \textstyle = \sigma_{0}^{2} + \frac{\sigma_{0} \theta_{0}(2 \tau \theta_{0} - \sin(2 \tau \theta_{0})) + \tau^{2} \theta_{0}^{4} + \theta_{0}^{2} \sin^{2}(\tau \theta_{0})}{\tau^{2} \theta_{0}^{2} - \sin^{2}(\tau \theta_{0} )}, \displaybreak[0] \label{CoeffsComplexA0} \\
\alpha_1 & \textstyle = 2 \theta_{0} e^{\sigma_{0} \tau} \frac{\tau \theta_{0} \cos(\tau \theta_{0} ) - \sin(\tau \theta_{0} )}{\tau^{2} \theta_{0}^{2} - \sin^{2}(\tau \theta_{0} )}, \displaybreak[0] \label{CoeffsComplexAlpha1} \\
\alpha_0 & \textstyle = 2 \theta_{0} e^{\sigma_{0} \tau} \frac{(\sigma_{0} - \tau \theta_0^2) \sin(\tau \theta_{0} ) - \tau \sigma_{0} \theta_{0} \cos(\tau \theta_{0} )}{\tau^{2} \theta_{0}^{2} - \sin^{2}(\tau \theta_{0} )} \label{CoeffsComplexAlpha0}
\end{align}
\end{subequations}
Moreover, if \eqref{CoeffsComplex} is satisfied, then all roots $\lambda$ of $\Delta$ which are different from $\lambda_0$ and from $\overline \lambda_0$ satisfy $\Re(\lambda) < \Re(\lambda_0)$.
\end{theorem}

As $\theta_0 \to 0$, the expressions in \eqref{CoeffsComplex} converge to those for the existence of a real root of multiplicity $4$ at $\lambda = \sigma_0$, as stated in Theorem~\ref{thm:MID}, for which it is known that $\sigma_0$ is dominant. The proof of Theorem~\ref{thm:mid-complex}, presented in \cite{MBNV:20}, uses this fact and the continuity of roots as functions of $\theta_0$ in order to show that dominance of the roots $\lambda_0$ and $\overline\lambda_0$ is preserved as $\theta_0$ increases from $0$ to any positive value, taking appropriate care of possible roots coming from $\infty$.
\begin{figure}[!ht]
\centering
\includegraphics[width=0.42\columnwidth]{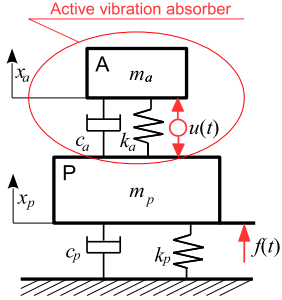}
\caption{Primary structure (P), with an active vibration absorber (A) to suppress displacement $x_p$ induced by harmonic disturbance force $f(t)$}
\label{fig:intro_1}
\end{figure}
Based on \cite{MBNV:20}, we illustrate the application of Theorem~\ref{thm:mid-complex} to \emph{active vibration suppression (AVS)}, in which we desire to suppress vibrations from an excitation force of a known frequency $\omega$ (see Fig.~\ref{fig:intro_1}). The system main body is a vibrating platform $P$ excited by a periodic external force $f(t) = F \cos(\omega t)$, and the absorber $A$ is actuated with the active feedback $u(t)$ to compensate the vibrations. The absorber dynamics is then $\ddot x_a(t)+2\zeta\Omega \dot x_a(t) + \Omega^2 x_a(t) = \frac{1}{m_a} u(t)$.

\begin{remark}
The \emph{delayed resonator scheme} by \cite{olgac1994novel} consists in guaranteeing vibration suppression by ensuring the overall system to have zeros at $\pm j \omega$, and we adapt that scheme here by placing zeros of multiplicity two at $\pm j \omega$ using Theorem~\ref{thm:mid-complex}. For a different adaptation, we refer to \cite{Kure2018Spectral}. $\Box$
\end{remark}

Choosing the delay as $\tau_k = \frac{k\pi}{\omega}$ for some $k \in \mathbb N^\ast$, Theorem~\ref{thm:mid-complex} ensures that $\pm j \omega$ is a \emph{double root\/} of $\Delta$ if and only if
\begin{equation}
\label{Delta_w}
\Delta_\omega(\lambda) = \lambda^2 + \frac{2}{\tau_k}\left((-1)^k e^{- \lambda \tau_k}-1\right)s + \omega^2.
\end{equation}
Introducing the active feedback $u(t)=m_a(\Omega^2-\omega^2)x_a(t)+2m_a\bigl(\zeta\Omega+\frac{1}{\tau_k}\bigr)\dot x_a(t) -2m_a\frac{(-1)^k}{\tau_k}\dot x_a(t-\tau_k)$, the characteristic function of the active absorber is given by \eqref{Delta_w}, with a double root at $\pm j\omega$. As shown, e.g., in \cite{Kure2018Spectral}, the transfer function $f\rightarrow x_p$ is in the form $G(\lambda)=\frac{\Delta_\omega(\lambda)}{M(\lambda)}$. Therefore, as required, the double roots at $\pm j\omega$ become double zeros of $G_{x_af}$, implying that no vibrations are transferred from $f$ to $x_p$ and the platform is fully silenced.

\section{Concluding remarks}
\label{sec:notes}

This paper presents several approaches and methods for handling multiple characteristic roots in time-delay systems represented by linear DDEs. More precisely, perturbation theory techniques, frequency-sweeping based approach and multiplicity-induced-dominancy method are explicitly discussed. For a better understanding of the concepts, notions as well as of the proposed methods, several illustrative examples complete the presentation. 

\section*{Acknowledgements}
\label{sec:acks}
The authors wish thank {\sc Jie Chen, Keqin Gu, Tamas Insperger, Wim Michiels, Rifat Sipahi} and  {\sc Tomas Vyhlidal} for useful discussions concerning the topics of the paper. C.F. {\sc M\'{e}ndez-Barr\'{i}os}  gratefully acknowledges the financial support of CentraleSup\'{e}lec (France) and CONACyT (Mexico) which supported several short stays of the author in France (2018-2021). {\sc X.-G. Li}'s work was supported in part by Fundamental
Research Funds for the Central Universities (No.
N2104007).

\bibliography{BIB-TDS-2021}         

\end{document}